\documentclass[11pt,reqno]{amsart}
   \usepackage[centertags]{amsmath}
   \usepackage{amsfonts}
   \usepackage{amssymb,mathtools}
   \usepackage{amsthm}
   \usepackage[utf8]{inputenc}
   \usepackage{enumitem,xcolor}
\usepackage{array}

\usepackage{caption}
\usepackage{graphicx}
\usepackage{subcaption}

\usepackage{hyperref}

\usepackage{geometry}

\usepackage{tikz}

\usepackage{color}

\usepackage{bm}

\newcommand\D{\displaystyle}

\allowdisplaybreaks[3]

\theoremstyle{plain}
\newtheorem{theorem}{Theorem}[section]

\theoremstyle{definition}

\theoremstyle{remark}
\newtheorem{remark}[theorem]{Remark}
\newtheorem*{remark*}{Remark}

\numberwithin{equation}{section}

\title[]{Darboux transformations of linear birth--death processes}

\author{Jos\'e Arcia-Manoleskos}
\address{Jos\'e Arcia-Manoleskos, Departamento de An\'{a}lisis Matem\'{a}tico. Universidad de Sevilla. Apdo (P. O. BOX) 1160. 41080 Sevilla. Spain.}
\address{Departamento de Matem\'{a}tica. Facultad de Ciencias Naturales, Exactas y Tecnolog\'{i}a. Universidad de Panam\'a}
\email{josarcman@alum.us.es}
\email{jose.arciam@up.ac.pa}

\author{Manuel Dom\'inguez de la Iglesia}
\address{Manuel Dom\'inguez de la Iglesia, Departamento de Matem\'aticas, Universidad de C\'adiz, Spain.}
\email{manuel.dominguezdelaiglesia@uca.es}

\thanks{This work was partially supported by research project PID2024-155133NB-I00 funded by MICIU/AEI, Spain, and PIUAH25/CC-006, funded by the Universidad de Alcal\'a, Spain. The first author was supported by a doctoral scholarship from SENACYT, Panama, under program BIDP-2022-026.}

\date{\today}

\makeatletter
   \@namedef{subjclassname@2020}{\textup{2020} Mathematics Subject Classification}
   \makeatother



\begin{document}

\begin{abstract}
In a previous paper we characterized the lower--upper (LU) and upper--lower (UL) factorizations of the infinitesimal generator of a possibly nonconservative birth--death process whose Darboux transformations are again generators of birth--death processes. Here we apply these results to several classical models with linear rates and explicitly known spectral measures: the $M/M/\infty$ queue (Charlier polynomials), the linear birth--death processes (Meixner and Laguerre polynomials), and the linear models of Ismail, Letessier and Valent (associated Laguerre polynomials). For each of them we compute the factors in closed form, determine the admissible ranges of the parameters, and describe the rates, spectral measures and orthogonal polynomials of the Darboux-transformed processes. 
\end{abstract}

\maketitle

\section{Introduction}

Birth--death processes are continuous-time Markov chains on $\mathbb{N}_0$ in which transitions occur only between neighboring states. Since the seminal work of Karlin and McGregor \cite{KMc2,KMc3}, it is well known that these processes admit a complete spectral description in terms of a spectral measure and an associated family of orthogonal polynomials, which yields an integral representation of the transition probabilities and spectral characterizations of recurrence, absorption and many other probabilistic properties (see \cite{MDIB} and also \cite{Scho}).

\smallskip

LU-type factorizations and the corresponding Darboux transformations, obtained by reversing the order of the factors, have been extensively used in the theory of orthogonal polynomials and integrable systems (see, for instance, \cite{GrH1,GrH3,GrHH}), and, in the discrete-time setting, as a systematic tool for constructing new birth--death chains with closely related spectral measures \cite{GdI}. In \cite{AdI} we extended this approach to continuous-time birth--death processes, possibly with absorption at a state $-1$ accessible from $0$. More precisely, given the infinitesimal generator $\mathcal{A}$ of such a process, we considered lower--upper (LU) and upper--lower (UL) factorizations of $\mathcal{A}$ and characterized those for which the corresponding Darboux transformations are again infinitesimal generators of birth--death processes. In the LU case, the factorization is uniquely determined once the absorption rate $\hat\mu_0$ of the transformed process is fixed, whereas the UL factorization involves an additional free parameter $x_0$. We obtained closed-form expressions for all the coefficients of the factors, determined the exact admissible ranges of the parameters, gave probabilistic interpretations of the coefficients, and showed that the spectral measures of the Darboux-transformed processes are given by Christoffel (LU case) and Geronimus (UL case) transformations of the original spectral measure. The general theory was illustrated in \cite{AdI} with the $M/M/1$ queue, whose constant birth and death rates lead to a spectral measure and orthogonal polynomials of Chebyshev type.

\smallskip

The purpose of the present paper is to study in detail two further classical families of birth--death processes with linear rates, namely the $M/M/\infty$ queue \cite{KMc4} and the linear birth--death processes \cite{KMc5,ILV}, including the models of \cite{ILV} related to the associated Laguerre polynomials. Their spectral measures are given, respectively, by the Poisson distribution and, depending on the relation between the birth and death parameters, by negative binomial or Gamma distributions, and their orthogonal polynomials are the Charlier, Meixner and Laguerre polynomials (and the associated Laguerre polynomials for the models of \cite{ILV}). For these examples, all the quantities involved in the factorizations can be computed explicitly, and several interesting phenomena appear:

\begin{itemize}[leftmargin=0.25in]
\item In the LU case, the coefficients of the factors and the rates of the Darboux-transformed processes can be written in closed form in terms of incomplete Gamma and Beta functions or Gauss hypergeometric functions. The Christoffel transformation shifts the parameters of the classical measures, and the corresponding Darboux birth--death polynomials are again Charlier, Meixner or Laguerre polynomials. In the positive recurrent cases, the two extreme admissible values of $\hat\mu_0$ lead to birth--death rates written in terms of the distribution function or its tail of the stationary distribution, and are related to Doob $h$-transforms of the dual processes.
\smallskip
\item In the UL case, the factorization cannot be applied directly to the recurrent processes with $\mu_0=0$, but it can be applied to their dual processes. The resulting spectral measures are affine combinations of the original spectral measure and a Dirac mass at the origin. The admissible range of the absorption rate
$\tilde\mu_0$ is governed by a series $T$.
For the models considered here, we obtain explicit
representations of this series in terms of the
exponential integral, generalized hypergeometric
functions, or definite integrals. In several
particular cases, further simplifications are
possible, including representations in terms of
the dilogarithm. We also investigate the convergence,
asymptotic behavior and monotonicity properties
of $T$ with respect to the model parameters. 
\end{itemize}
In all the examples considered here, none of the original processes or of their Darboux transforms is explosive, since their birth rates grow at most linearly, and all the Darboux-transformed processes again have a uniquely determined spectral measure (see Section~\ref{sec2}).

\smallskip

The structure of the paper is as follows. In Section~\ref{sec2} we collect the notation and the results of \cite{AdI} that will be used throughout the paper. Section~\ref{secMMinf} is devoted to the $M/M/\infty$ queue, and Section~\ref{seclinear} to the linear birth--death processes, where the three cases $\lambda<\mu$ (Meixner), $\lambda>\mu$ (Meixner) and $\lambda=\mu$ (Laguerre) are treated separately. In Section~\ref{secassLag} we study the two linear models of \cite{ILV}, related to the associated Laguerre polynomials, which differ only in the absorption rate at $0$; in this case the LU and UL factorizations connect both models, and they also show that the spectral measure of the conservative model has a mass at $0$ in the positive recurrent case. Finally, Section~\ref{secConc2} contains some concluding remarks and open questions.

\section{Preliminaries}\label{sec2}

In this section we collect the notation and the main results of \cite{AdI} that will be used in the rest of the paper. We refer to \cite{AdI} for the proofs and further details.

\subsection{Birth--death processes and their spectral representation}

Let $\{X_t:t\geq0\}$ be a birth--death process on $\mathbb{N}_0$, possibly killed at the absorbing state $-1$, with birth and death rates $\{\lambda_n,\mu_n\}$ and infinitesimal generator
\begin{equation}\label{QQmm}
\mathcal{A}=
	\begin{pmatrix}
		-(\lambda_0+\mu_0) & \lambda_0 &  &  & \\
		\mu_1  & -(\lambda_1+\mu_1)  & \lambda_1 &  &\\
		 & \mu_2 & -(\lambda_2+\mu_2) & \lambda_2 & \\
		&& \ddots & \ddots & \ddots
	\end{pmatrix}.
\end{equation}
$\mathcal{A}$ is conservative if and only if $\mu_0=0$; when $\mu_0>0$, the process may jump from state $0$ to the absorbing state $-1$. We denote by $\mathcal{A}^{(\alpha)}$ the generator obtained from \(\mathcal A\) by replacing \(\mu_0\) with \(\alpha\), leaving all other rates unchanged. The \emph{potential coefficients} are given by
\begin{equation}\label{potcoef}
\pi_0=1,\quad\pi_n=\frac{\lambda_0\lambda_1\cdots\lambda_{n-1}}{\mu_1\mu_2\cdots\mu_n},\quad n\geq1,
\end{equation}
and we will use the quantities
\begin{equation}\label{quantity_AB}
A=\sum_{n=0}^{\infty}\frac{1}{\lambda_n\pi_n},
\quad
B=\sum_{n=0}^{\infty}\pi_n,
\quad
S=\sum_{n=0}^{\infty}
\pi_n
\left(
1+\mu_0\sum_{k=0}^{n-1}\frac{1}{\lambda_k\pi_k}
\right)^2.
\end{equation}
As in \cite{AdI}, we assume that
\begin{equation}\label{uniqueness_assumption}
\begin{cases}
A+B=\infty, & \mu_0=0,\\[1mm]
S=\infty, & \mu_0>0,
\end{cases}
\end{equation}
which is equivalent to the determinacy of the associated Stieltjes moment problem (see \cite{KMc2,vanDoorn2003,vanDoorn2016}), and $\{X_t:t\geq0\}$ denotes the minimal process associated with $\mathcal{A}$. Then there exists a unique probability measure $d\psi$ supported on $[0,\infty)$ and a family of polynomials $(Q_n)_{n\in\mathbb{N}_0}$ defined by
\begin{equation}\label{QQs}
Q_{-1}(x)=0,\quad Q_0(x)=1,\quad
-xQ_n(x)
=
\lambda_nQ_{n+1}(x)
-(\lambda_n+\mu_n)Q_n(x)
+\mu_nQ_{n-1}(x),
\quad n\geq0,
\end{equation}
such that the transition probabilities $P_{ij}(t)=\mathbb{P}(X_t=j\mid X_0=i)$ are given by the \emph{Karlin--McGregor formula}
\begin{equation*}
P_{ij}(t)
=
\pi_j\int_0^\infty
e^{-xt}Q_i(x)Q_j(x)\,d\psi(x),
\quad i,j\in\mathbb{N}_0.
\end{equation*}
Evaluating \eqref{QQs} at $x=0$ we have
\begin{equation}\label{Qn0}
Q_n(0) = 1 + \mu_0 \sum_{k=0}^{n-1} \frac{1}{\lambda_k \pi_k}, 
\quad n \ge 0,
\end{equation}
so that $S=\sum_{n\ge0}\pi_nQ_n^2(0)$. We will also use the $0$-th associated polynomials $(Q_n^{(0)})_{n\in\mathbb{N}_0}$, defined by the recurrence \eqref{QQs} for $n\ge1$ with initial conditions $Q_0^{(0)}(x)=0$ and $Q_1^{(0)}(x)=-1/\lambda_0$, for which
\begin{equation}\label{Qn00}
Q_n^{(0)}(0) = -\sum_{k=0}^{n-1} \frac{1}{\lambda_k \pi_k}, 
\quad n \ge 1.
\end{equation}
Recall that, if $\mu_0=0$, the process is recurrent if and only if $A=\infty$, and positive recurrent if and only if, in addition, $B<\infty$; if $\mu_0>0$, absorption at $-1$ is certain if and only if $A=\infty$ (see \cite{KMc3,MDIB}). Finally, given a strictly positive sequence $h$ with $(\mathcal{A}h)(i)=0$ for all $i\in\mathbb{N}_0$, the \emph{Doob $h$-transform} of $\mathcal{A}$ (see \cite[Section 17.6]{LP}) is the operator
\begin{equation}\label{doobt}
(\mathcal{A}_h f)(i)=\frac{1}{h(i)}\mathcal{A}(h f)(i),
\quad i \in \mathbb{N}_0,
\end{equation}
which is the generator of a conservative birth--death process with rates $\lambda_i^h=\lambda_ih(i+1)/h(i)$, $i\ge0$, and $\mu_i^h=\mu_ih(i-1)/h(i)$, $i\ge1$.

\subsection{LU factorization}\label{subsecLU}

Following \cite[Section~3]{AdI}, consider the LU factorization $\mathcal{A}=\widetilde{\mathcal{A}}_L\widetilde{\mathcal{A}}_U$, where
\begin{equation*}
\widetilde{\mathcal{A}}_L=
		\begin{pmatrix}
			\tilde{s}_0 &  &  &  \\
			\tilde{r}_1 & \tilde{s}_1 &  &  \\
			 & \tilde{r}_2 &\tilde{s}_2 &  \\
			 & & \ddots    & \ddots
		\end{pmatrix},\quad
		\widetilde{\mathcal{A}}_U=
		\begin{pmatrix}
			\tilde{y}_0 & \tilde{x}_0 &  &  & \\
			 & \tilde{y}_1 & \tilde{x}_1 &  & \\
			&  &\tilde{y}_2 &  \tilde{x}_2 & \\
			 &  &  & \ddots & \ddots
		\end{pmatrix},
\end{equation*}
normalized so that $\tilde s_0=1$, and such that the Darboux transformation $\widehat{\mathcal{A}}=\widetilde{\mathcal{A}}_U\widetilde{\mathcal{A}}_L$ satisfies $\widehat{\mathcal{A}}\bm e=-\hat\mu_0e^{(0)}$ for some $\hat\mu_0\ge0$, where $\bm e$ is the vector of all ones and $e^{(0)}=(1,0,0,\ldots)^T$. Once $\hat\mu_0$ is fixed, the coefficients are uniquely determined. More precisely, writing $\tilde q_n=\tilde y_n/\tilde x_n$, $n\ge0$, and $\tilde t_n=\tilde s_n+\tilde r_n$, $n\ge1$, $\tilde t_0=1$, we have (see \cite[Proposition 3.1]{AdI})
\begin{equation*}
\tilde q_n= -\frac{Q_{n+1}(0)}{Q_n(0)},\quad n\geq0,\quad
\tilde t_n = \left(1-\frac{\hat\mu_0}{\lambda_0+\mu_0}\right)Q_n(0),\quad n\geq1,
\end{equation*}
and
\begin{align}
\label{ssnn}\tilde s_n&=\frac{1}{\pi_n Q_n(0)}\left[\frac{\hat\mu_0}{\lambda_0+\mu_0}+\left(1-\frac{\hat\mu_0}{\lambda_0+\mu_0}\right)\sum_{k=0}^n \pi_kQ_k^2(0)\right], \quad n\ge1,\quad \tilde s_0=1,\\
\label{rrnn}\tilde r_n&=-\frac{1}{\pi_n Q_n(0)}\left[\frac{\hat\mu_0}{\lambda_0+\mu_0}+\left(1-\frac{\hat\mu_0}{\lambda_0+\mu_0}\right)\sum_{k=0}^{n-1} \pi_kQ_k^2(0)\right], \quad n\ge1,\\
\label{xxnn}\tilde x_n&=\frac{\lambda_n\pi_nQ_n(0)}{\dfrac{\hat\mu_0}{\lambda_0+\mu_0}+\left(1-\dfrac{\hat\mu_0}{\lambda_0+\mu_0}\right)\D\sum_{k=0}^n \pi_kQ_k^2(0)},\quad n\ge1,\quad \tilde x_0=\lambda_0,\\
\label{yynn}\tilde y_n&=-\tilde{x}_n\frac{Q_{n+1}(0)}{Q_n(0)}=-\frac{\lambda_n\pi_nQ_{n+1}(0)}{\dfrac{\hat\mu_0}{\lambda_0+\mu_0}+\left(1-\dfrac{\hat\mu_0}{\lambda_0+\mu_0}\right)\D\sum_{k=0}^n \pi_kQ_k^2(0)},\quad n\ge1,\quad \tilde y_0=-(\lambda_0+\mu_0).
\end{align}
The birth--death rates of the Darboux transformation $\widehat{\mathcal{A}}$ are given by
\begin{equation}\label{newbdcoe}
\hat\lambda_n=\lambda_n\frac{\tilde s_{n+1}}{\tilde s_n},\quad\hat\mu_{n+1}=\frac{\lambda_{n+1}\mu_{n+1}}{\lambda_n}\cdot\frac{\tilde q_{n+1}}{\tilde q_n}\cdot\frac{\tilde s_n}{\tilde s_{n+1}},\quad n\geq0.
\end{equation}
The following result is one of the main results of \cite{AdI} (see \cite[Theorem 3.6]{AdI}).

\begin{theorem}\label{theom33}
Assume the standing condition \eqref{uniqueness_assumption}, and consider the LU factorization $\mathcal{A}=\widetilde{\mathcal{A}}_L\widetilde{\mathcal{A}}_U$ with coefficients \eqref{ssnn} and \eqref{rrnn}--\eqref{yynn}. Then $\widehat{\mathcal A}=\widetilde{\mathcal{A}}_U \widetilde{\mathcal{A}}_L$ is the generator of a birth--death process with birth--death rates \eqref{newbdcoe} if and only if the absorption rate $\hat\mu_0$ of the Darboux-transformed process satisfies
\begin{enumerate}
\item if $\mu_0>0$:
\begin{equation*}
    0 \le \hat\mu_0 \le \lambda_0+\mu_0;
\end{equation*}
\item if $\mu_0=0$:
\begin{equation}\label{condLU2}
0\le\hat\mu_0\le
\begin{cases}
\dfrac{\lambda_0B}{B-1}, & B<\infty,\\[2mm]
\lambda_0, & B=\infty.
\end{cases}
\end{equation}
\end{enumerate}
\end{theorem}
The birth--death polynomials $(\hat Q_n)_{n\in\mathbb{N}_0}$ associated with $\widehat{\mathcal A}$ are given by
\begin{equation}\label{Qnhat}
\hat{Q}_n(x) = -\frac{1}{x}\left(\tilde{x}_n Q_{n+1}(x) + \tilde{y}_n Q_n(x)\right), \quad n \ge 0,
\end{equation}
and they are orthogonal with respect to the \emph{Christoffel transform} of $d\psi$ (see \cite{AdI} and also \cite{GdI,Yo})
\begin{equation}\label{spmec1}
d\widehat{\psi}(x)=\frac{x}{\lambda_0+\mu_0}\,d\psi(x).
\end{equation}

\begin{remark}\label{remdual}
If $\mu_0=0$ and $\hat\mu_0=\lambda_0$, the coefficients of the LU factorization reduce to $\tilde s_n=1/\pi_n$, $\tilde r_{n+1}=-1/\pi_{n+1}$, $\tilde x_n=\lambda_n\pi_n$ and $\tilde y_n=-\lambda_n\pi_n$, $n\ge0$, and the birth--death rates of the Darboux transformation $\widehat{\mathcal{A}}=\widehat{\mathcal{A}}^{(\lambda_0)}$ are $\hat\lambda_n=\mu_{n+1}$, $\hat\mu_n=\lambda_n$, $n\ge0$. This is the \emph{dual process} of \cite{KMc2}, whose generator will be denoted by $\mathcal{A}^d$ (see \cite[Remark 3.15]{MDIB}).
\end{remark}

Assume now that $\mu_0=0$ and $B<\infty$, so that the process is positive recurrent, and let $X_\pi$ be a random variable with the stationary distribution $\mathbb{P}(X_\pi=n)=\pi_n/B$, $n\ge0$, and cumulative distribution function $F_{X_\pi}(n)=\mathbb{P}(X_\pi\le n)$. Then \eqref{ssnn} can be written as (see \cite[Section~3.1]{AdI})
\begin{equation}\label{ssnnB}
\tilde{s}_n
=
\frac{1}{\mathbb{P}(X_\pi=n)}
\left[
\frac{\hat{\mu}_0}{\lambda_0B}
+
\left(1-\frac{\hat{\mu}_0}{\lambda_0}\right)
F_{X_\pi}(n)
\right],
\quad n\geq0,
\end{equation}
and, if $\hat\mu_0=0$, the birth--death rates of $\widehat{\mathcal{A}}$ take the form
\begin{equation}\label{adfsdf}
\hat{\lambda}_n
=
\mu_{n+1}\frac{F_{X_\pi}(n+1)}{F_{X_\pi}(n)},
\quad
\hat{\mu}_{n+1}
=
\lambda_{n+1}\frac{F_{X_\pi}(n)}{F_{X_\pi}(n+1)},
\quad n\geq0,
\end{equation}
while, if $\hat\mu_0=\lambda_0B/(B-1)$, the same formulas hold after replacing $F_{X_\pi}$ by the tail function $1-F_{X_\pi}$. The sequence $h(n)=F_{X_\pi}(n)$ is harmonic with respect to the dual process $\mathcal{A}^d$, so that $\widehat{\mathcal{A}}^{(0)}$ is precisely the Doob $h$-transform \eqref{doobt} of $\mathcal{A}^d$ with $h(n)=F_{X_\pi}(n)$. In contrast, the tail function $h(n)=1-F_{X_\pi}(n)$ is harmonic only in the interior, so that the case $\hat\mu_0=\lambda_0B/(B-1)$ should not be regarded as a global Doob $h$-transform of $\mathcal{A}^d$, although the corresponding $h$-transform relation remains valid for the interior rates (see \cite[Section~3.1]{AdI}).

\subsection{UL factorization}\label{subsecUL}

Following \cite[Section~4]{AdI}, consider now the UL factorization $\mathcal{A}=\mathcal{A}_U\mathcal{A}_L$, where
\begin{equation*}
		\mathcal{A}_U=
		\begin{pmatrix}
			y_0 & x_0 &  &  & \\
			 &y_1 & x_1 &  & \\
			&  &y_2 &  x_2 & \\
			 &  &  & \ddots & \ddots
		\end{pmatrix},\quad 
		\mathcal{A}_L=
		\begin{pmatrix}
			s_0 &  &  &  \\
			r_1 & s_1 &  &  \\
			 & r_2 &s_2 &  \\
			 & & \ddots    & \ddots
		\end{pmatrix},
\end{equation*}
normalized so that $s_0=1$, and such that the Darboux transformation $\widetilde{\mathcal{A}}=\mathcal{A}_L\mathcal{A}_U$ satisfies $\widetilde{\mathcal{A}}\bm e=-\tilde\mu_0e^{(0)}$ for some $\tilde\mu_0\ge0$. In this case the factorization depends on a \emph{free parameter} $x_0$. Setting
\begin{equation}\label{uun}
u_n=Q_n(0)+(x_0+\tilde\mu_0)Q_n^{(0)}(0),\quad n\geq0,
\end{equation}
the coefficients of the factors are given by (see \cite[Section~4]{AdI})
\begin{align}
\label{xxnn2}x_n&=\frac{1}{\pi_nu_n}\left[x_0+\tilde\mu_0 -\tilde\mu_0\sum_{k=0}^n \pi_ku_k^2\right], \quad n\ge0,\\
\label{yynn2}y_n&=-\frac{1}{\pi_nu_n}\left[x_0 +\tilde\mu_0-\tilde\mu_0\sum_{k=0}^{n-1} \pi_ku_k^2\right],\quad n\ge1,\quad y_0=-(x_0+\tilde\mu_0),\\
\label{ssnn2}s_n&=\frac{\lambda_{n-1}\pi_{n-1}u_{n-1}}{x_0 +\tilde\mu_0-\tilde\mu_0\displaystyle\sum_{k=0}^{n-1} \pi_ku_k^2},\quad n\ge1,\quad s_0=1,\\
\label{rrnn2}r_n&=-\frac{\mu_{n}\pi_{n}u_{n}}{x_0+\tilde\mu_0 -\tilde\mu_0\displaystyle\sum_{k=0}^{n-1} \pi_ku_k^2},\quad n\ge1,
\end{align}
and the birth--death rates of the Darboux transformation $\widetilde{\mathcal{A}}$ are given by
\begin{equation}\label{newbdcoe2}
\tilde\lambda_n
=
\lambda_{n-1}\frac{x_n}{x_{n-1}},
\quad
\tilde\mu_{n+1}
=
\frac{\lambda_n\mu_n}{\lambda_{n-1}}
\frac{u_{n+1}u_{n-1}}{u_n^2}
\frac{x_{n-1}}{x_n},\quad n\geq1,
\end{equation}
together with the boundary rates
\begin{equation}\label{laconcha3r}
\tilde\lambda_0=x_0,
\quad
\tilde\mu_1
=
\frac{(\lambda_0+\mu_0-x_0-\tilde\mu_0)(x_0+\tilde\mu_0)}{x_0}.
\end{equation}
The following result is the second main result of \cite{AdI} (see \cite[Theorem 4.5]{AdI}). Here and below we use the convention $1/\infty=0$.

\begin{theorem}\label{theom37}
Assume the standing condition \eqref{uniqueness_assumption}, consider the UL factorization $\mathcal{A}=\mathcal{A}_U\mathcal{A}_L$ with coefficients \eqref{xxnn2}--\eqref{rrnn2}, and define the series
\begin{equation}\label{TTT}
T=\sum_{n=0}^\infty\pi_nu_n^2=\sum_{n=0}^\infty\pi_n\left(1+(\mu_0-x_0-\tilde\mu_0)\sum_{k=0}^{n-1}\frac{1}{\lambda_k\pi_k}\right)^2.
\end{equation}
Then $\widetilde{\mathcal A}=\mathcal{A}_L \mathcal{A}_U$ is the generator of a birth--death process with birth and death rates given by \eqref{newbdcoe2} and \eqref{laconcha3r} if and only if
\begin{equation}\label{condULgen}
0<x_0+\tilde\mu_0\leq \mu_0+\frac{1}{A},
\quad
\tilde\mu_0\geq0,\quad\mbox{together with}\quad \begin{cases}
\tilde\mu_0 T\leq x_0+\tilde\mu_0, & T<\infty,\\[1mm]
\tilde\mu_0=0, & T=\infty.
\end{cases}
\end{equation}
Under these conditions, the admissible parameters split into the following two cases:
\begin{enumerate}
\item If $x_0+\tilde\mu_0\neq\mu_0+\frac{1}{A},$ then $T=\infty$. In this situation $\tilde\mu_0=0,$ and the free parameter $x_0$ is subject to the condition
\begin{equation}\label{cond1UL}
0<x_0<\mu_0+\frac{1}{A}.
\end{equation}
\item If $x_0+\tilde\mu_0=\mu_0+\frac{1}{A},$ then $T$ is given by
\begin{equation}\label{TTT3}
T=\begin{cases}
\displaystyle
\sum_{n=0}^\infty \pi_n
\left(1-\frac{1}{A}\sum_{k=0}^{n-1}\frac{1}{\lambda_k\pi_k}\right)^2,
& \mbox{if } A<\infty,\\[1ex]
B,
& \mbox{if } A=\infty,
\end{cases}
\end{equation}
and the absorption rate $\tilde\mu_0\geq0$ must satisfy
\begin{equation}\label{condBAB}
0 \le \tilde{\mu}_0 \le \frac{1}{T}\left(\mu_0+\frac{1}{A}\right).
\end{equation}
In the case where $T=\infty,$ then $\tilde\mu_0=0$ and $x_0=\mu_0+\frac{1}{A}$.
\end{enumerate}
\end{theorem}
The birth--death polynomials $(\tilde Q_n)_{n\in\mathbb{N}_0}$ associated with $\widetilde{\mathcal A}$ are given by
\begin{equation}\label{QQbb1}
\tilde Q_0(x)=1,\quad\tilde{Q}_n(x) = r_n Q_{n-1}(x) + s_n Q_n(x), \quad n \geq 1,
\end{equation}
or, equivalently (see \cite[Section~4]{AdI}),
\begin{equation}\label{QDTUL}
\tilde{Q}_n(x)=\frac{x_0+\tilde\mu_0-x\D\sum_{k=0}^{n-1}\pi_ku_kQ_k(x)}{x_0+\tilde\mu_0-\tilde\mu_0\D\sum_{k=0}^{n-1}\pi_ku_k^2},\quad n\geq1,
\end{equation}
and they are orthogonal with respect to the \emph{Geronimus transform} of $d\psi$ (see \cite{AdI} and also \cite{GdI,Yo})
\begin{equation}\label{spmec2e}
d\widetilde{\psi}(x)
=\frac{x_0+\tilde\mu_0}{x}\,d\psi(x)+M\,\delta_0(dx),
\quad
M=1-(x_0+\tilde\mu_0)m_{-1},
\end{equation}
where \(\delta_0\) denotes the Dirac measure at \(0\) and $m_{-1} = \int_{0}^{\infty} x^{-1} d\psi(x)$ is given by (see \cite[(3.70)]{MDIB})
\begin{equation}\label{mmm1}
m_{-1}=\frac{A}{1+\mu_0A}=\left(\mu_0+\frac{1}{A}\right)^{-1}.
\end{equation}
Condition \eqref{condULgen} guarantees that $M\ge0$. A probabilistic interpretation of the coefficients \eqref{xxnn2}--\eqref{rrnn2} in terms of hitting probabilities and expected passage times can be found in \cite[Section~4.1]{AdI}.

\begin{remark}\label{remdual2}
If $\mu_0>0$, $\tilde\mu_0=0$ and $x_0=\mu_0$, the coefficients of the UL factorization reduce to $s_{n+1}=\lambda_n\pi_n/\mu_0$, $r_{n+1}=-\lambda_n\pi_n/\mu_0$, $x_n=\mu_0/\pi_n$ and $y_n=-\mu_0/\pi_n$, $n\ge0$, and the birth--death rates of the Darboux transformation $\widetilde{\mathcal{A}}=\widetilde{\mathcal{A}}^{(0)}$ are $\tilde\lambda_n=\mu_n$, $\tilde\mu_{n+1}=\lambda_n$, $n\ge0$. This is the dual process in the case $\mu_0>0$ (see \cite[Remark 3.12]{MDIB}). In contrast, when $\mu_0=0$, the UL factorization of the dual generator $\mathcal A^d$, with $x_0=\lambda_0$ and $\tilde\mu_0=0$, recovers the original generator $\mathcal A$.
\end{remark}

\begin{remark}\label{remsimpl}
If $\tilde\mu_0=0$, then \eqref{xxnn2} gives $x_n=x_0/(\pi_nu_n)$, and \eqref{newbdcoe2}--\eqref{laconcha3r} reduce to
\[
\tilde\lambda_0=x_0,\quad \tilde\lambda_n=\mu_n\frac{u_{n-1}}{u_n},\ n\ge1,\quad \tilde\mu_{n+1}=\lambda_n\frac{u_{n+1}}{u_n},\ n\ge0.
\]
Similarly, if $\mu_0=0$, the LU rates \eqref{newbdcoe} are $\hat\lambda_n=\mu_{n+1}\tilde d_{n+1}/\tilde d_n$ and $\hat\mu_{n+1}=\lambda_{n+1}\tilde d_n/\tilde d_{n+1}$, with $\tilde d_n=\frac{\hat\mu_0}{\lambda_0}+\bigl(1-\frac{\hat\mu_0}{\lambda_0}\bigr)\sum_{k=0}^n\pi_k$, and $\hat\lambda_n\le\lambda_n+\mu_{n+1}$ for every admissible $\hat\mu_0$. Indeed, writing $\kappa=\hat\mu_0/\lambda_0$, we have $\tilde d_{n+1}=\tilde d_n+(1-\kappa)\pi_{n+1}$ with $\tilde d_n>0$. If $\kappa\le1$, then $\tilde d_n\ge(1-\kappa)\pi_n$, so that $\tilde d_{n+1}/\tilde d_n\le1+\pi_{n+1}/\pi_n$ and, since $\mu_{n+1}\pi_{n+1}=\lambda_n\pi_n$, $\hat\lambda_n\le\mu_{n+1}+\lambda_n$. If $\kappa>1$, then $(\tilde d_n)_{n\in\mathbb{N}_0}$ is decreasing and $\hat\lambda_n\le\mu_{n+1}$.
\end{remark}

\subsection{Determinacy and nonexplosion of the Darboux-transformed processes}

By \cite[Proposition~4.8]{AdI}, the Darboux transform $\widetilde{\mathcal{A}}$ obtained from the UL factorization always satisfies the standing condition \eqref{uniqueness_assumption}. For the LU factorization, $\widehat{\mathcal{A}}$ satisfies \eqref{uniqueness_assumption} if and only if
\begin{equation}\label{sww}
\widehat S=(\lambda_0+\mu_0)\sum_{n=0}^{\infty}\frac{\Sigma_n^2}{\lambda_n\pi_nQ_n(0)Q_{n+1}(0)}=\infty,\quad \Sigma_n=\sum_{k=0}^n\pi_kQ_k^2(0),
\end{equation}
a condition which does not depend on $\hat\mu_0$ (see \cite[Proposition~3.8]{AdI}). In the LU factorizations of Sections~\ref{secMMinf} and~\ref{seclinear}, and of Model II in Section~\ref{secassLag}, we have $\mu_0=0$, so that $Q_n(0)=1$ and $\Sigma_n\ge\max\{1,\pi_n\}$, and therefore $\widehat S\ge\lambda_0\sum_n1/\lambda_n=\infty$, since the birth rates grow at most linearly. For Model I in Section~\ref{secassLag}, where $\mu_0>0$, a similar argument applies. Hence all the Darboux-transformed processes considered below again have a uniquely determined spectral measure. Moreover, a birth--death process with $\sum_n1/\lambda_n=\infty$ is not explosive, since Reuter's series $\sum_n(\lambda_n\pi_n)^{-1}\sum_{k\le n}\pi_k$ is bounded below by $\sum_n1/\lambda_n$ (for more details, see \cite{And}). The same criterion applies when \(\mu_0>0\) with an absorbing state \(-1\). For the LU factorizations with $\mu_0=0$ the divergence of \(\sum_n 1/\hat{\lambda}_n\) follows from the bound $\hat\lambda_n\le\lambda_n+\mu_{n+1}$ in Remark~\ref{remsimpl}, and for the remaining factorizations, including the UL factorizations of the dual processes in Sections~\ref{secMMinf} and~\ref{seclinear}, from the explicit rates given below, which grow at most linearly. Hence none of these processes is explosive, and the minimal process is the only birth--death process with the given rates.

\section{The $M/M/\infty$ queue}\label{secMMinf}

Let $\{X_t, t\geq0\}$ be the birth--death process with constant birth and linear death rates of the form
$$
\lambda_n=\lambda,\quad n\geq0,\quad\mu_n=n\mu,\quad n\geq0,\quad\lambda,\mu>0.
$$
The potential coefficients are given by
\begin{equation*}
    \pi_n=\frac{(\lambda/\mu)^n}{n!}, \quad n\geq 0.
\end{equation*}
It is well known (see \cite[Section 3]{KMc4} or \cite[Example 3.36]{MDIB}) that the spectral measure associated with the infinitesimal generator $\mathcal{A}$ is given by
\[
d\psi(x)=e^{-\lambda/\mu}\sum_{n=0}^{\infty}\frac{(\lambda/\mu)^{n}}{n!}\delta_{\mu n}(dx),
\]
where \(\delta_a\) denotes the Dirac measure at \(a\). Observe that the support of the measure is $\{0,\mu,2\mu,\ldots\}$ and 
the masses are the Poisson probabilities with parameter $\lambda/\mu$. The corresponding orthogonal polynomials can be written in terms of the \emph{Charlier polynomials} $C_n(x;a),n\geq0$. Indeed,
\[
Q_n(x)=C_n\left(\frac{x}{\mu};\frac{\lambda}{\mu}\right),\quad n\ge0.
\]
Since the process is conservative, i.e. $\mu_0=0$, then \(Q_n(0)=1\). One readily checks that the quantities $A$ and $B$ in \eqref{quantity_AB} are given by $A=\infty$ and $B=e^{\lambda/\mu}$. 
Therefore, the process is positive recurrent. Let $Y_\pi$ be a random variable with the stationary distribution
\[
\mathbb{P}(Y_\pi=n)
=
e^{-\lambda/\mu}\frac{(\lambda/\mu)^n}{n!},
\quad n\geq0.
\]
Thus, $Y_\pi\sim\operatorname{Poisson}(\lambda/\mu)$. In this case,
\begin{equation*}
\sum_{k=0}^n\pi_k
=
e^{\lambda/\mu}\mathbb{P}(Y_\pi\leq n)
=
e^{\lambda/\mu}\frac{\Gamma(n+1,\lambda/\mu)}{n!},
\end{equation*}
where $\Gamma(s,x)=\int_x^\infty t^{s-1}e^{-t}\,dt$ denotes the \emph{upper incomplete Gamma function}. We shall write
\(
F_{Y_\pi}(n)=\mathbb{P}(Y_\pi\leq n)
\)
for the cumulative distribution function of $Y_\pi$.

\medskip

\hspace{-.35cm}\underline{\textbf{LU factorization}}. The coefficients $\tilde s_n, \tilde r_n, \tilde x_n,$ and $\tilde y_n$ in \eqref{ssnn} and \eqref{rrnn}--\eqref{yynn} can be explicitly computed. Given that \(\mu_0=0\), we have that $\tilde q_n=-1$ and $\tilde t_n=1-\hat\mu_0/\lambda, \tilde t_0=1$. Therefore
\begin{equation}\label{srxymmi}
\begin{split}
\tilde s_n&=\frac{\hat\mu_0\mu^nn!}{\lambda^{n+1}}+\left(1-\frac{\hat\mu_0}{\lambda}\right)\frac{\mu^n}{\lambda^n}e^{\lambda/\mu}\Gamma(n+1,\lambda/\mu),\quad n\geq0,\\
\tilde r_n&=-\frac{\hat\mu_0\mu^nn!}{\lambda^{n+1}}-\left(1-\frac{\hat\mu_0}{\lambda}\right)\frac{n\mu^n}{\lambda^n}e^{\lambda/\mu}\Gamma(n,\lambda/\mu),\quad n\geq1,\\
\tilde x_n&=\frac{\lambda}{\tilde s_n}=-\tilde y_n,\quad n\geq0.
\end{split}
\end{equation}
These coefficients can also be expressed in terms of the stationary random variable $Y_\pi$ and its cumulative distribution function $F_{Y_\pi}(n)$, as in \eqref{ssnnB}. In this setting, the coefficient $\tilde{s}_n$ can be written as
\[
\tilde{s}_n
=
\frac{1}{\mathbb{P}(Y_\pi=n)}
\left[
F_{Y_\pi}(n)
+
\frac{\hat{\mu}_0}{\lambda}
\bigl(e^{-\lambda/\mu}-F_{Y_\pi}(n)\bigr)
\right],
\quad n\geq0.
\]
With this, we can compute the coefficients $\hat\lambda_n,\hat\mu_{n+1},n\geq0,$ in \eqref{newbdcoe}, given by
$$
\hat\lambda_n=\lambda\frac{\tilde s_{n+1}}{\tilde s_n},\quad\hat\mu_{n+1}=\mu(n+1)\frac{\tilde s_n}{\tilde s_{n+1}},\quad n\geq0.
$$
In order to obtain the Darboux transformation \(\widehat{\mathcal{A}}\) as the infinitesimal generator of a new birth--death process, we apply the second part of Theorem \ref{theom33}, since $\mu_0=0$. Moreover, since $B=e^{\lambda/\mu}$, the parameter \(\hat\mu_0\) must satisfy condition~\eqref{condLU2}, namely,
\[
0\le\hat\mu_0\le\frac{\lambda}{1-e^{-\lambda/\mu}}.
\]
As for the Darboux spectral measure, since \(\tilde y_0=-\lambda\), we have
\begin{equation}\label{smdbs}
d\widehat\psi(x)=\frac{x}{\lambda}d\psi(x)=e^{-\lambda/\mu}\sum_{n=0}^\infty\frac{\mu n(\lambda/\mu)^n}{\lambda n!}\delta_{\mu n}(dx)=e^{-\lambda/\mu}\sum_{n=0}^\infty\frac{(\lambda/\mu)^n}{n!}\delta_{\mu (n+1)}(dx),
\end{equation}
that is, the transformed measure has the original Poisson weights shifted by one lattice step, so that its support is \(\{\mu,2\mu,3\mu,\ldots\}\), with the point \(0\) removed from the original support. Using that $\tilde y_n=-\tilde x_n$, \eqref{Qnhat}, \cite[(9.14.6)]{KS} and the self-duality of Charlier polynomials, we obtain that the Darboux birth--death polynomials \((\hat Q_n)_{n\in\mathbb{N}_0}\) are given, as expected, by
    \begin{equation}\label{darbchar}
        \hat Q_n(x)=\frac{\tilde x_n}{\lambda}C_n\left(\frac{x}{\mu}-1;\frac{\lambda}{\mu}\right),\quad n\geq0.
        \end{equation}
We now examine some particular cases where further simplifications arise.
\begin{itemize}[leftmargin=0.25in]
\item Case $\hat\mu_0=0$. In this case it follows, from \eqref{srxymmi}, that
$$
\tilde s_n=\frac{\mu^ne^{\lambda/\mu}}{\lambda^n}\Gamma(n+1,\lambda/\mu),\quad \tilde r_n=-\frac{n\mu^ne^{\lambda/\mu}}{\lambda^n}\Gamma(n,\lambda/\mu),\quad\tilde x_n=\frac{\lambda^{n+1}e^{-\lambda/\mu}}{\mu^n\Gamma(n+1,\lambda/\mu)}=-\tilde y_n,
$$
where $\tilde s_n,\tilde x_n,\tilde y_n$ are defined for $n\geq0$, whereas $\tilde r_n$ is defined for $n\geq1$. Therefore
\begin{equation*}
\hat\lambda_n=\mu\frac{\Gamma(n+2,\lambda/\mu)}{\Gamma(n+1,\lambda/\mu)},\quad \hat\mu_{n+1}=\lambda(n+1)\frac{\Gamma(n+1,\lambda/\mu)}{\Gamma(n+2,\lambda/\mu)},\quad n\geq0.
\end{equation*}
These birth--death rates can also be written (see \eqref{adfsdf}) as
\[
\hat{\lambda}_n
 = \mu(n+1)\frac{F_{Y_\pi}(n+1)}{F_{Y_\pi}(n)},
 \quad
\hat{\mu}_{n+1}
 = \lambda\frac{F_{Y_\pi}(n)}{F_{Y_\pi}(n+1)},
 \quad n\geq 0,
\]
where $F_{Y_\pi}(n)$ denotes the cumulative distribution function of the stationary random variable $Y_\pi$.

\item Case $\hat\mu_0=\lambda(1-e^{-\lambda/\mu})^{-1}$. This limiting case is also interesting since, from \eqref{srxymmi} and after some computations, we obtain
$$
\tilde s_n=\frac{\mu^n\gamma(n+1,\lambda/\mu)}{(1-e^{-\lambda/\mu})\lambda^n},\quad \tilde r_n=-\frac{n\mu^n\gamma(n,\lambda/\mu)}{(1-e^{-\lambda/\mu})\lambda^n},\quad\tilde x_n=\frac{\lambda^{n+1}(1-e^{-\lambda/\mu})}{\mu^n\gamma(n+1,\lambda/\mu)}=-\tilde y_n,
$$
where $\tilde s_n,\tilde x_n,\tilde y_n$ are defined for $n\geq0$, whereas $\tilde r_n$ is defined for $n\geq1$. Here $\gamma(s,x)=\int_0^x t^{s-1}e^{-t}dt$ is the \emph{lower incomplete Gamma function}. Therefore
\begin{equation*}
\hat\lambda_n=\mu\frac{\gamma(n+2,\lambda/\mu)}{\gamma(n+1,\lambda/\mu)},\quad \hat\mu_{n+1}=\lambda(n+1)\frac{\gamma(n+1,\lambda/\mu)}{\gamma(n+2,\lambda/\mu)},\quad n\geq0.
\end{equation*}
In this case, these birth--death rates can also be written in terms of the tail distribution function of the stationary random variable $Y_\pi$.
\end{itemize}
As mentioned after \eqref{adfsdf}, the operator $\widehat{\mathcal{A}}$ corresponding to $\hat{\mu}_0=0$ (that is, $\widehat{\mathcal{A}}^{(0)}$) coincides with the \emph{Doob $h$-transform} (see~\eqref{doobt}) of the dual process $\mathcal{A}^d$ (or $\widehat{\mathcal{A}}^{(\lambda)}$), obtained by taking $h(n)=F_{Y_\pi}(n)$. The birth--death rates of the dual process are given by $\hat{\lambda}_n=\mu(n+1)$ and $\hat{\mu}_n=\lambda$ for all $n\geq0$. The other limiting case, corresponding to $\hat{\mu}_0=\lambda(1-e^{-\lambda/\mu})^{-1}$, is related instead to the tail function $h(n)=1-F_{Y_\pi}(n)$. In this case, $h$ is harmonic only in the interior, and hence  \(\widehat{\mathcal{A}}\) is not a global Doob $h$-transform of $\mathcal{A}^d$ (see the end of Section~\ref{subsecLU}). 

\medskip

\hspace{-.35cm}\underline{\textbf{UL factorization}}. Since \(A=\infty\) and \(\mu_0=0\), an admissible UL factorization of the type considered here is not possible, since condition \eqref{condULgen} with $\mu_0=0$ and $A=\infty$ would force $x_0+\tilde\mu_0\le0$, which is is not allowed. However, the factorization becomes possible if we instead apply the UL decomposition to the dual process \(\mathcal{A}^d=\widehat{\mathcal{A}}^{(\lambda)}\), whose birth--death rates are
\(\hat{\lambda}_n=\mu(n+1),\hat{\mu}_{n}=\lambda,n\geq 0\). This process is nonconservative ($\mu_0=\lambda$) and can be interpreted as a single-server queue with service rate $\lambda$ whose arrival rate grows linearly with the number of customers present, killed at rate $\lambda$ when the system is empty. Abusing notation, we denote the dual birth--death rates \(\{\hat{\lambda}_n,\hat{\mu}_n\}\) by \(\{\lambda_n,\mu_n\}\), so that \(\lambda_n=\mu(n+1)\) and \(\mu_n=\lambda,\) for \(n\geq 0\). The corresponding potential coefficients are given by
$$
\pi_n=n!\left(\frac{\mu}{\lambda}\right)^n,\quad n\geq 0.
$$
A straightforward computation shows that the quantities \(A\) and \(B\) in \eqref{quantity_AB} are given by
\[
A=\frac{1}{\lambda}\left(e^{\lambda/\mu}-1\right),\quad B=\infty.
\]
Since \(A<\infty\), absorption at \(-1\) is not certain. Moreover, a direct computation using \eqref{Qn0} (see also \eqref{darbchar} evaluated at \(x=0\) and \(\tilde{x}_n=\lambda(\lambda/\mu)^n/n!\)) yields
$$
Q_n(0)=e^{\lambda/\mu}\dfrac{\Gamma(n+1,\lambda/\mu)}{n!},\quad n\ge0,
$$
where \(\Gamma(s,x)\) denotes the upper incomplete Gamma function. Since $Q_n(0)\to e^{\lambda/\mu}>0$ as $n\to\infty$ and $B=\infty$, it follows that the quantity \(S\) in \eqref{quantity_AB} diverges. Consequently, the associated moment problem is determinate.

\begin{figure}[htbp]
\centering
\begin{subfigure}[t]{0.48\textwidth}
    \centering
    \includegraphics[width=\textwidth]{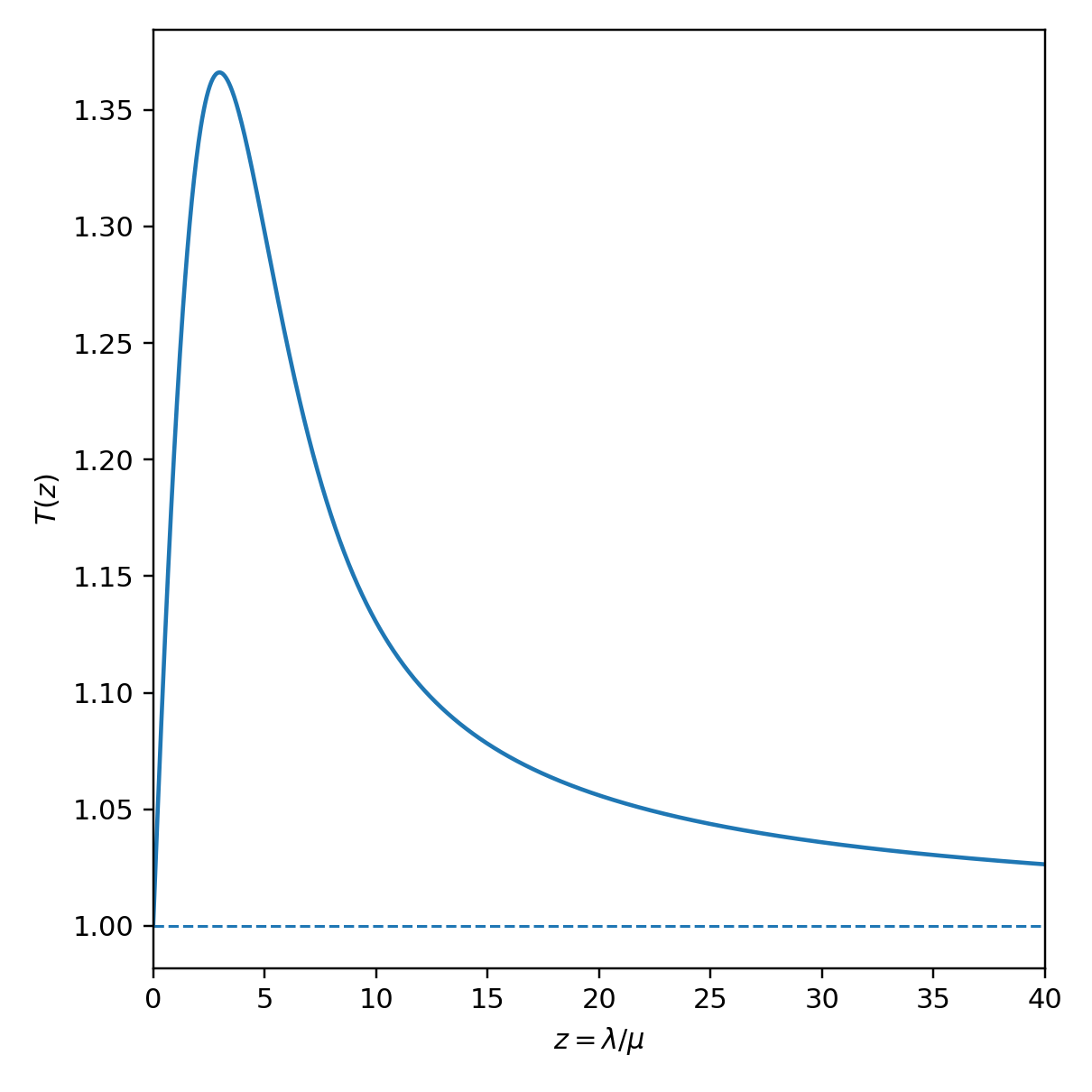}\end{subfigure}
\hfill
\begin{subfigure}[t]{0.48\textwidth}
    \centering
    \includegraphics[width=\textwidth]{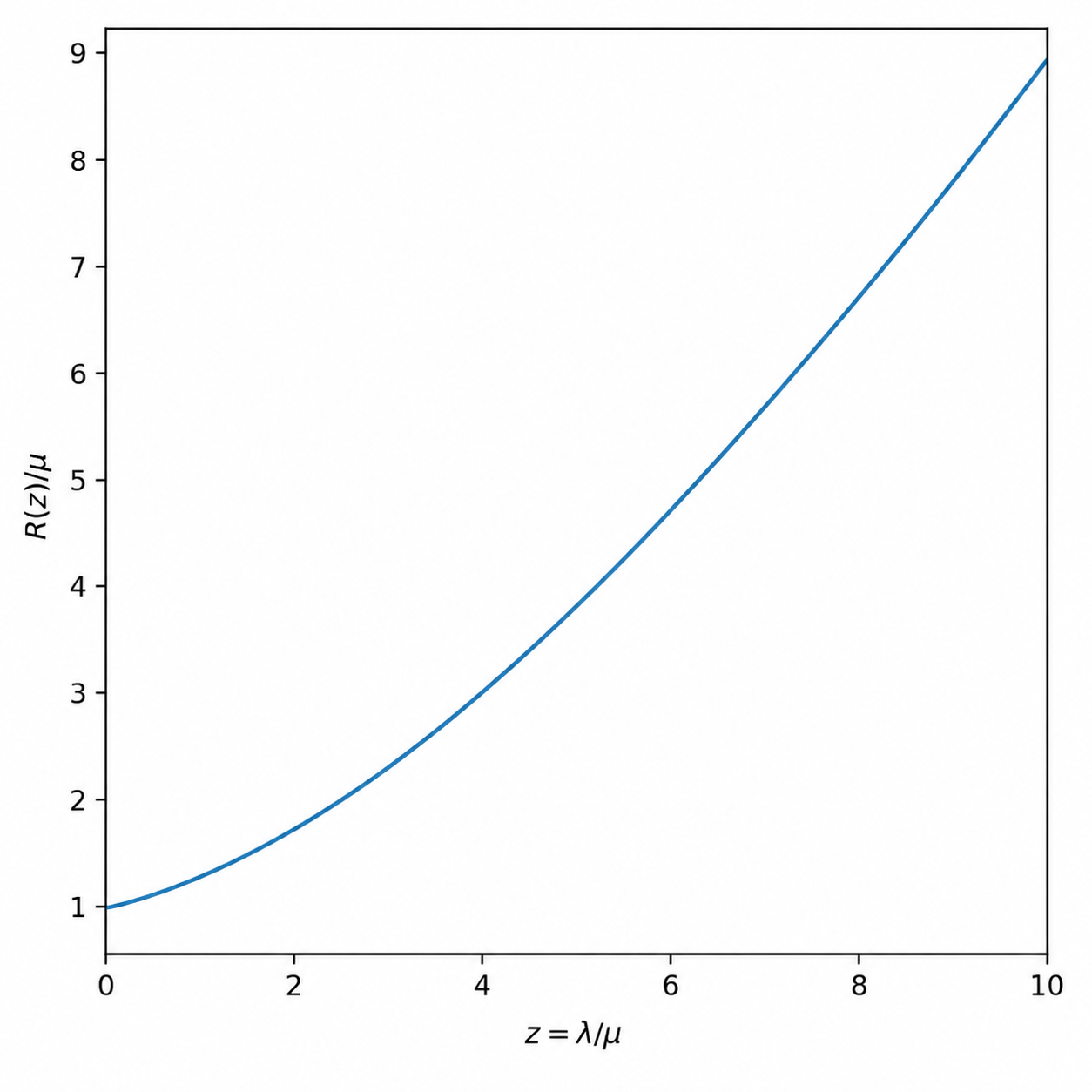}
\end{subfigure}
\caption{Graphs of $T(z)$ and $R(z)/\mu$ as functions of $z=\lambda/\mu$.}
\label{fig:T_function}
\end{figure}

Let us now perform the UL factorization of this dual process. By Theorem \ref{theom37}, if \(x_0+\tilde{\mu}_0\neq\frac{\lambda e^{\lambda/\mu}}{e^{\lambda/\mu}-1}\), then \(T=\infty\), and consequently \(\tilde{\mu}_0=0\). In this case, there is a \emph{free parameter} \(x_0\), which must satisfy
\[
0 < x_0 < \frac{\lambda e^{\lambda/\mu}}{e^{\lambda/\mu}-1}.
\]
Here $u_n=\bigl(1-\frac{x_0}{\lambda}\bigr)e^{\lambda/\mu}\frac{\Gamma(n+1,\lambda/\mu)}{n!}+\frac{x_0}{\lambda}$, and by Remark~\ref{remsimpl} the rates are $\tilde\lambda_0=x_0$, $\tilde\lambda_n=\lambda\,u_{n-1}/u_n$ and $\tilde\mu_{n+1}=\mu(n+1)\,u_{n+1}/u_n$. On the other hand, if \(x_0+\tilde{\mu}_0=\frac{\lambda e^{\lambda/\mu}}{e^{\lambda/\mu}-1}\), then the series \(T\) in \eqref{TTT3} (corresponding to the case \(A<\infty\)) can be expressed as
$$
T=\left(\frac{e^{\lambda/\mu}}{e^{\lambda/\mu}-1}\right)^2 \sum_{n=0}^{\infty}\frac{(\mu/\lambda)^n}{n!}\gamma(n+1,\lambda/\mu)^2,
$$
where \(\gamma(s,x)\) denotes the lower incomplete Gamma function. Using the integral representation of \(\gamma(s,x)\) and performing some straightforward computations, we obtain
\begin{equation*}
T=
\frac{\lambda}{\mu}
\frac{e^{\lambda/\mu}}{\left(e^{\lambda/\mu}-1\right)^2}
\left[
\mathrm{Ei}\!\left(\frac{\lambda}{\mu}\right)
-\gamma_{\mathrm{E}}
-\log\!\left(\frac{\lambda}{\mu}\right)
\right],
\end{equation*}
where 
\[
\mathrm{Ei}(x)
=
\gamma_{\mathrm{E}}
+\log x
+\int_0^x \frac{e^t-1}{t}\,dt,\quad x>0,
\]
is the \emph{exponential integral} and \(\gamma_{\mathrm{E}}\) is the Euler--Mascheroni constant. Therefore, according to \eqref{condBAB}, the absorption rate \(\tilde{\mu}_0\) must satisfy
\begin{equation*}
0\leq\tilde\mu_0 \le\frac{\mu(e^{\lambda/\mu}-1)}{\operatorname{Ei}\left(\frac{\lambda}{\mu}\right)-\gamma_{\mathrm{E}}-\log\left(\frac{\lambda}{\mu}\right)}.
\end{equation*}
If we call $R=\frac{1}{T}\left(\lambda+\frac{1}{A}\right)$, the dependence of \(T=T(z)\) and $R/\mu=R(z)/\mu$ on the ratio \(z=\lambda/\mu\) are illustrated in Figure \ref{fig:T_function}.

Since $\mathrm{Ei}(z)-\gamma_{\mathrm E}-\log z\sim e^z/z$ as $z\to\infty$, we have $T(z)=1+z/4+\mathcal O(z^2)$ as $z\to0^+$ and $T(z)\to1$ as $z\to\infty$; in particular $T$ is not monotone, and it attains its maximum $T\approx1.366$ at $z\approx2.97$.

As for the Darboux spectral measure, the original spectral measure of the dual process is given by \eqref{smdbs}. We must therefore consider the Geronimus transformation of \eqref{smdbs}, namely \eqref{spmec2e}. From \eqref{mmm1} it follows that
$$
m_{-1}=\frac{1}{\lambda}\left(1-e^{-\lambda/\mu}\right).
$$
Substituting this expression into \eqref{spmec2e} and rearranging the resulting series, we obtain after straightforward simplifications
\[
d\widetilde\psi(x)
=
\frac{x_0+\tilde\mu_0}{\lambda}e^{-\lambda/\mu}
\sum_{n=1}^{\infty}\frac{(\lambda/\mu)^n}{n!}\delta_{\mu n}(dx)
+
\left[
1-\frac{x_0+\tilde\mu_0}{\lambda}
\left(1-e^{-\lambda/\mu}\right)
\right]\delta_0(dx).
\]
Equivalently,
\[
d\widetilde\psi(x)
=
\frac{x_0+\tilde\mu_0}{\lambda}\,d\psi(x)
+
\left(1-\frac{x_0+\tilde\mu_0}{\lambda}\right)\delta_0(dx),
\]
where $d\psi$ is the original spectral measure associated with the Charlier polynomials. Thus, the masses of $d\psi$ at the positive points $\mu,2\mu,\ldots$ are rescaled by the factor $(x_0+\tilde\mu_0)/\lambda$, while the remaining mass is concentrated at $x=0$. As for the orthogonal polynomials, a lengthy but straightforward computation, using \cite[(9.14.6)]{KS} and the self-duality of Charlier polynomials together with the identity
$\Gamma(n+1,\lambda/\mu)=n\Gamma(n,\lambda/\mu)+(\lambda/\mu)^ne^{-\lambda/\mu}$,
yields $\tilde Q_0=1$ and
\begin{align*}
\tilde Q_n(x)
&=
\frac{s_n}{n!}\left(\frac{\lambda}{\mu}\right)^n
\left[
C_n\left(\frac{x}{\mu};\frac{\lambda}{\mu}\right)
-
\left(1-\frac{x_0+\tilde\mu_0}{\lambda}\right)
\frac{(\lambda/\mu)^{n-1}}{(n-1)!u_{n-1}}
C_{n-1}\left(\frac{x}{\mu}-1;\frac{\lambda}{\mu}\right)
\right],
\quad n\geq1.
\end{align*}
In particular, when $x_0+\tilde\mu_0=\lambda$, we recover the original Charlier spectral measure and the corresponding classical Charlier polynomials.

\section{Linear birth--death processes}\label{seclinear}

Let  $\{X_t, t\geq0\}$ be the birth--death process with linear birth--death rates 
$$
\lambda_n=(n+\beta)\lambda,\quad n\geq0,\quad\mu_n=n\mu,\quad n\geq0,\quad\lambda,\mu,\beta>0.
$$
The potential coefficients are given by
\begin{equation*}
    \pi_n=\frac{(\beta)_n(\lambda/\mu)^n}{n!}, \quad n\geq 0,
\end{equation*}
where $(a)_n=a(a+1)\cdots(a+n-1), (a)_0=1,$ denotes the usual Pochhammer symbol. Observe that the infinitesimal generator \(\mathcal{A}\) is always conservative, that is, \(\mu_0 = 0\). Therefore $Q_n(0)=1$. The spectral analysis associated with \(\mathcal{A}\) can be found in \cite{KMc5} (see also \cite[Example~3.38]{MDIB}). More general linear models, including nonconservative ones, were studied in \cite{ILV}, where the associated Laguerre and Meixner polynomials appear; the associated Laguerre case will be treated below in Section~\ref{secassLag}.

The quantities $A$ and $B$ in~\eqref{quantity_AB} are given by
\begin{equation}\label{A_B_linear}
\begin{split}
A&=\frac{1}{\lambda\beta}\sum_{n=0}^{\infty}\frac{n!}{(\beta+1)_n}\left(\frac{\mu}{\lambda}\right)^n=\frac{1}{\lambda\beta}\,{}_{2}F_{1}\!\left(\!\begin{array}{c}
1,1 \\[-1.5mm]
\beta+1
\end{array}; \frac{\mu}{\lambda} \right),\\
B&=\sum_{n=0}^{\infty}\frac{(\beta)_n}{n!}\left(\frac{\lambda}{\mu}\right)^n={}_{1}F_{0}\!\left(\!\begin{array}{c}
\beta\\[-1.5mm]
-
\end{array}; \frac{\lambda}{\mu} \right),
\end{split}
\end{equation}
where ${_p}F_q$ is the \emph{generalized hypergeometric function} (also denoted by \({_p}F_q(a_1,\hdots,a_p;b_1,\hdots,b_q;x)\)). The hypergeometric notation above is understood in terms of the corresponding defining series whenever it converges; when the series diverges, the corresponding quantity \(A\) or \(B\) is \(+\infty\). Depending on the values of \(\lambda\) and \(\mu\), the nature of the series \(A\) and \(B\) may vary, in some situations also depending on the parameter \(\beta\). We distinguish three cases: \(\lambda < \mu\), \(\lambda > \mu\), and \(\lambda = \mu\). Each case will be treated separately, since both the spectral measure and the associated orthogonal polynomials differ substantially among them. 
Moreover, the probabilistic behavior of the corresponding birth--death process also changes significantly across these three regimes.

\subsection{Case $\lambda<\mu$.} Since $\lambda/\mu<1$ and $\frac{n!}{(\beta+1)_n}\sim\Gamma(\beta+1)n^{-\beta}$ as $n\to\infty$, the terms of \(A\) in \eqref{A_B_linear} grow geometrically, so $A=\infty$ and $B$ is the binomial series, that is,
$$
B=\left(1-\frac{\lambda}{\mu}\right)^{-\beta}.
$$
Therefore, in this situation, the birth--death process is positive recurrent. It is well known that the spectral measure associated with the infinitesimal generator $\mathcal{A}$ is given by
\begin{equation*}
    d\psi(x)=\left(1-\frac{\lambda}{\mu}\right)^{\beta}\sum_{n=0}^{\infty}\frac{(\beta)_n (\lambda/\mu)^n}{n!}\delta_{(\mu-\lambda)n}(dx),
\end{equation*}
where \(\delta_a\) denotes the Dirac measure at \(a\). Observe that the support of the measure is $\{0,\mu-\lambda,2(\mu-\lambda),\ldots\}$ and the masses are the negative binomial or Pascal probabilities with parameters $r=\beta$ and $p=1-\lambda/\mu$. The corresponding orthogonal polynomials can be written in terms of the \emph{Meixner polynomials} $M_n(x;b,c),n\geq0$. Indeed,
\begin{equation*}
    Q_n(x)=M_n\left(\frac{x}{\mu-\lambda};\beta,\frac{\lambda}{\mu}\right),\quad n\geq0.
\end{equation*}
Since the process is positive recurrent, let $Z_\pi$ be a random variable with the stationary distribution
\[
\mathbb{P}(Z_\pi=n)
=
\left(1-\frac{\lambda}{\mu}\right)^{\beta}
\frac{(\beta)_n}{n!}
\left(\frac{\lambda}{\mu}\right)^n,
\quad n\geq0.
\]
Observe that $Z_\pi\sim\mbox{NB}(\beta,1-\lambda/\mu)$ (that is, a negative binomial random variable). In this case, we also have that
\begin{equation}\label{posss2}
\sum_{k=0}^n\pi_k=\left(1-\frac{\lambda}{\mu}\right)^{-\beta}\mathbb{P}(Z_\pi\leq n)=\left(1-\frac{\lambda}{\mu}\right)^{-\beta}\frac{\mbox{B}(1-\lambda/\mu;\beta,n+1)}{\mbox{B}(\beta,n+1)},
\end{equation}
where $\mbox{B}(x;a,b)=\int_0^xt^{a-1}(1-t)^{b-1}dt$ is the \emph{incomplete Beta function} and $\mbox{B}(a,b)=\mbox{B}(1;a,b)$ is the \emph{Beta function}. Observe that $F_{Z_\pi}(n)=\mathbb{P}(Z_\pi\leq n)$ is the cumulative distribution function of the negative binomial random variable $Z_\pi$.

\medskip

\hspace{-.35cm}\underline{\textbf{LU factorization}}. The coefficients $\tilde s_n, \tilde r_n, \tilde x_n,$ and $\tilde y_n$ in \eqref{ssnn} and \eqref{rrnn}--\eqref{yynn} can be explicitly computed. Given that \(\mu_0=0\), we have that $\tilde q_n=-1$ and $\tilde t_n=1-\frac{\hat\mu_0}{\lambda\beta}, \tilde t_0=1$. Therefore, using that $\mbox{B}(\beta,n+1)=n!/(\beta)_{n+1}$, we have
\begin{equation}\label{srxymmi2}
\begin{split}
\tilde s_n&=\frac{\hat\mu_0\mu^nn!}{\beta(\beta)_n\lambda^{n+1}}+\left(1-\frac{\hat\mu_0}{\lambda\beta}\right)\frac{(n+\beta)\mu^n}{\lambda^n}\left(1-\frac{\lambda}{\mu}\right)^{-\beta}\mbox{B}(1-\lambda/\mu;\beta,n+1),\quad n\geq0,\\
\tilde r_n&=-\frac{\hat\mu_0\mu^nn!}{\beta(\beta)_n\lambda^{n+1}}-\left(1-\frac{\hat\mu_0}{\lambda\beta}\right)\frac{n\mu^n}{\lambda^n}\left(1-\frac{\lambda}{\mu}\right)^{-\beta}\mbox{B}(1-\lambda/\mu;\beta,n),\quad n\geq1,\\
\tilde x_n&=\frac{\lambda(n+\beta)}{\tilde s_n}=-\tilde y_n,\quad n\geq0.
\end{split}
\end{equation}
These coefficients can also be expressed in terms of the stationary random variable $Z_\pi$ and its cumulative distribution function $F_{Z_\pi}(n)$, as in \eqref{ssnnB}. In this setting, the coefficient $\tilde{s}_n$ can be written as
\[
\tilde{s}_n
=
\frac{1}{\mathbb{P}(Z_\pi=n)}
\left[
F_{Z_\pi}(n)
+
\frac{\hat{\mu}_0}{\lambda\beta}
\left(
\left(1-\frac{\lambda}{\mu}\right)^{\beta}
-
F_{Z_\pi}(n)
\right)
\right],
\quad n\geq0.
\]
With this, we can compute the coefficients $\hat\lambda_n,\hat\mu_{n+1},n\geq0,$ in \eqref{newbdcoe}, given by
$$
\hat\lambda_n=\lambda(n+\beta)\frac{\tilde s_{n+1}}{\tilde s_n},\quad\hat\mu_{n+1}=\frac{\mu(n+\beta+1)(n+1)}{n+\beta}\frac{\tilde s_n}{\tilde s_{n+1}},\quad n\geq0.
$$
In order to obtain the Darboux transformation \(\widehat{\mathcal{A}}\) as the infinitesimal generator of a new birth--death process, we apply the second part of Theorem \ref{theom33}, since $\mu_0=0$. Since $B=\left(1-\frac{\lambda}{\mu}\right)^{-\beta}$, the parameter \(\hat\mu_0\) must satisfy condition~\eqref{condLU2}, that is,
\[
0\le\hat\mu_0\le\frac{\lambda\beta}{1-\left(1-\dfrac{\lambda}{\mu}\right)^{\beta}}.
\]
As for the Darboux spectral measure, since \(\tilde y_0=-\lambda\beta\), we have
\begin{equation}\label{meixasdf}
\begin{split}
d\widehat\psi(x)&=\frac{x}{\lambda\beta}d\psi(x)=\left(1-\frac{\lambda}{\mu}\right)^{\beta}\sum_{n=0}^\infty\frac{(\mu-\lambda)n(\beta)_n(\lambda/\mu)^n}{\lambda\beta n!}\delta_{(\mu-\lambda)n}(dx)\\
&=\left(1-\frac{\lambda}{\mu}\right)^{\beta+1}\sum_{n=0}^{\infty}\frac{(\beta+1)_n (\lambda/\mu)^n}{n!}\delta_{(\mu-\lambda)(n+1)}(dx),
\end{split}
\end{equation}
that is, the transformed measure has the same support as \(d\psi\), except for the removal of the point \(0\), while its weights correspond to those of a negative binomial distribution with parameter \(\beta+1\), shifted by one lattice step. Using that $\tilde y_n=-\tilde x_n$, \eqref{Qnhat}, \cite[(9.10.6)]{KS} and the self-duality of Meixner polynomials, we obtain that the Darboux birth--death polynomials \((\hat Q_n)_{n\in\mathbb{N}_0}\) are given, as expected, by
    \begin{equation}\label{darbmeix}
        \hat Q_n(x)=\frac{\tilde x_n}{\lambda\beta}M_n\left(\frac{x}{\mu-\lambda}-1;\beta+1,\frac{\lambda}{\mu}\right),\quad n\geq0.
        \end{equation}
We now examine some particular cases where further simplifications arise.
\begin{itemize}[leftmargin=0.25in]
\item Case $\hat\mu_0=0$. In this case it follows, from \eqref{srxymmi2}, that
\begin{align*}
\tilde s_n&=\frac{(n+\beta)\mu^n}{\lambda^n}\left(1-\frac{\lambda}{\mu}\right)^{-\beta}\mbox{B}(1-\lambda/\mu;\beta,n+1),\quad n\geq0,\\
\tilde r_n&=-\frac{n\mu^n}{\lambda^n}\left(1-\frac{\lambda}{\mu}\right)^{-\beta}\mbox{B}(1-\lambda/\mu;\beta,n),\quad n\geq1,\\
\tilde x_n&=\frac{\lambda^{n+1}}{\mu^n\mbox{B}(1-\lambda/\mu;\beta,n+1)}\left(1-\frac{\lambda}{\mu}\right)^{\beta}=-\tilde y_n,\quad n\geq0,
\end{align*}
and therefore
\begin{equation*}
\hat\lambda_n=\mu(n+\beta+1)\frac{\mbox{B}(1-\lambda/\mu;\beta,n+2)}{\mbox{B}(1-\lambda/\mu;\beta,n+1)},\quad \hat\mu_{n+1}=\lambda(n+1)\frac{\mbox{B}(1-\lambda/\mu;\beta,n+1)}{\mbox{B}(1-\lambda/\mu;\beta,n+2)},\quad n\geq0.
\end{equation*}
These birth--death rates can also be written (see \eqref{adfsdf}) as
\[
\hat{\lambda}_n
 = \mu(n+1)\frac{F_{Z_\pi}(n+1)}{F_{Z_\pi}(n)},
 \quad
\hat{\mu}_{n+1}
 = \lambda(n+\beta+1)\frac{F_{Z_\pi}(n)}{F_{Z_\pi}(n+1)},
 \quad n\geq 0,
\]
where $F_{Z_\pi}(n)$ denotes the cumulative distribution function of the 
stationary random variable $Z_\pi$.

\item Case $\hat\mu_0=\lambda\beta(1-(1-\lambda/\mu)^\beta)^{-1}$. This limiting case is also interesting since, from \eqref{srxymmi2} and after some computations, we obtain
\begin{align*}
\tilde s_n&=\frac{(n+\beta)\mu^n}{\lambda^n(1-(1-\lambda/\mu)^\beta)}\mbox{C}(1-\lambda/\mu;\beta,n+1),\quad n\geq0,\\
\tilde r_n&=-\frac{n\mu^n}{\lambda^n(1-(1-\lambda/\mu)^\beta)}\mbox{C}(1-\lambda/\mu;\beta,n),\quad n\geq1,\\
\tilde x_n&=\frac{\lambda^{n+1}(1-(1-\lambda/\mu)^\beta)}{\mu^n\mbox{C}(1-\lambda/\mu;\beta,n+1)}=-\tilde y_n,\quad n\geq0,
\end{align*}
where $\mbox{C}(x;a,b)=\mbox{B}(a,b)-\mbox{B}(x;a,b)=\int_x^1t^{a-1}(1-t)^{b-1}dt$ is the \emph{complementary incomplete Beta function}. Therefore
\begin{equation*}
\hat\lambda_n=\mu(n+\beta+1)\frac{\mbox{C}(1-\lambda/\mu;\beta,n+2)}{\mbox{C}(1-\lambda/\mu;\beta,n+1)},\quad \hat\mu_{n+1}=\lambda(n+1)\frac{\mbox{C}(1-\lambda/\mu;\beta,n+1)}{\mbox{C}(1-\lambda/\mu;\beta,n+2)},\quad n\geq0.
\end{equation*}
In this case, these birth--death rates can also be written in terms of the tail distribution function of the stationary random variable $Z_\pi$.
\end{itemize}
As mentioned after \eqref{adfsdf}, the operator $\widehat{\mathcal{A}}$ corresponding to $\hat{\mu}_0=0$ (that is, $\widehat{\mathcal{A}}^{(0)}$) coincides with the \emph{Doob $h$-transform} (see~\eqref{doobt}) of the dual process $\mathcal{A}^d$ (or $\widehat{\mathcal{A}}^{(\lambda\beta)}$), obtained by taking $h(n)=F_{Z_\pi}(n)$. The birth--death rates of the dual process are given by $\hat{\lambda}_n=\mu(n+1)$ and $\hat{\mu}_n=\lambda(n+\beta)$ for all $n\geq0$. The other limiting case, corresponding to $\hat\mu_0=\lambda\beta\bigl(1-(1-\lambda/\mu)^\beta\bigr)^{-1}$, is related instead to the tail function $h(n)=1-F_{Z_\pi}(n)$. In this case, $h$ is harmonic only in the interior, so that, again, this case \(\widehat{\mathcal{A}}\) is not a global Doob $h$-transform of $\mathcal{A}^d$ (see the end of Section~\ref{subsecLU}).

\medskip

\hspace{-.35cm}\underline{\textbf{UL factorization}}. Since \(A=\infty\) and \(\mu_0=0\), an admissible UL factorization of the type considered here is not possible, since condition \eqref{condULgen} with $\mu_0=0$ and $A=\infty$ would force $x_0+\tilde\mu_0\le0$, which is is not allowed. However, the factorization becomes possible if we instead apply the UL decomposition to the dual process \(\mathcal{A}^d=\widehat{\mathcal{A}}^{(\lambda\beta)}\), whose birth--death rates are \(\hat{\lambda}_n=\mu(n+1),\hat{\mu}_{n}=\lambda(n+\beta),n\geq 0\). Abusing notation, we denote the dual birth--death rates \(\{\hat{\lambda}_n,\hat{\mu}_n\}\) again by \(\{\lambda_n,\mu_n\}\), so that \(\lambda_n=\mu(n+1)\) and \(\mu_n=\lambda(n+\beta)\) for \(n\geq 0\). The corresponding potential coefficients are given by
$$
\pi_n=\frac{n!}{(\beta+1)_n}\left(\frac{\mu}{\lambda}\right)^n,\quad n\geq 0.
$$
A straightforward computation shows that the quantities \(A\) and \(B\) in \eqref{quantity_AB} are given by
\[
A=\frac{1}{\lambda\beta}\left[\left(1-\frac{\lambda}{\mu}\right)^{-\beta}-1\right],\quad B=\infty.
\]
Since \(A<\infty\), absorption at \(-1\) is not certain. Moreover, a direct computation using \eqref{Qn0} (see also \eqref{darbmeix} evaluated at \(x=0\) and \(\tilde{x}_n=\lambda(\beta)_{n+1}(\lambda/\mu)^n/n!\)) yields
\begin{equation}\label{Qn0mm}
Q_n(0)=\left(1-\frac{\lambda}{\mu}\right)^{-\beta}-\frac{(\beta)_{n+1}}{(n+1)!}\left(\frac{\lambda}{\mu}\right)^{n+1}{}_{2}F_{1}\!\left(\!\begin{array}{c}
1,\beta+n+1 \\[-1.5mm]
n+2
\end{array}; \frac{\lambda}{\mu} \right),\quad n\geq0.
\end{equation}
Since \(Q_n(0)\to\left(1-\frac{\lambda}{\mu}\right)^{-\beta}>0\) as \(n\to\infty\), it follows immediately that \(S\) in \eqref{quantity_AB} diverges. Therefore, the associated moment problem is determinate.

Let us now perform the UL factorization of this dual process. First of all, observe that 
$$
\mu_0+\frac{1}{A}=\frac{\lambda\beta}{1-\left(1-\frac{\lambda}{\mu}\right)^\beta}.
$$
By Theorem \ref{theom37}, if \(x_0+\tilde{\mu}_0\neq\frac{\lambda\beta}{1-\left(1-\frac{\lambda}{\mu}\right)^\beta}\), then \(T=\infty\), and consequently \(\tilde{\mu}_0=0\). In this case, there is a \emph{free parameter} \(x_0\), which must satisfy
\[
0 < x_0 < \frac{\lambda\beta}{1-\left(1-\frac{\lambda}{\mu}\right)^\beta}.
\]
Here $u_n=\bigl(1-\frac{x_0}{\lambda\beta}\bigr)Q_n(0)+\frac{x_0}{\lambda\beta}$, with $Q_n(0)$ given by \eqref{Qn0mm}, and by Remark~\ref{remsimpl} the rates are $\tilde\lambda_0=x_0$, $\tilde\lambda_n=\lambda(n+\beta)\,u_{n-1}/u_n$ and $\tilde\mu_{n+1}=\mu(n+1)\,u_{n+1}/u_n$. On the other hand, if \(x_0+\tilde{\mu}_0=\frac{\lambda\beta}{1-\left(1-\frac{\lambda}{\mu}\right)^\beta}\), then the series \(T\) in \eqref{TTT3} (corresponding to the case \(A<\infty\)) can be written as
$$
T=\frac{1}{{}_{2}F_{1}\!\left(\!\begin{array}{c}
1,\beta+1 \\[-1.5mm]
2
\end{array}; \D\frac{\lambda}{\mu} \right)^2}\sum_{n=0}^{\infty}\frac{(\beta+1)_{n}(\lambda/\mu)^n}{n!(n+1)^2}\,{}_{2}F_{1}\!\left(\!\begin{array}{c}
1,\beta+n+1 \\[-1.5mm]
n+2
\end{array}; \D\frac{\lambda}{\mu} \right)^2.
$$
The series $T$ actually converges, since $\frac{(\beta+1)_n}{n!}\sim\frac{n^{\beta}}{\Gamma(\beta+1)}$ and the part with the hypergeometric function ${}_{2}F_{1}$ is bounded as $n\to\infty$. Therefore, the sequence inside the sum behaves like $Cn^{\beta-2}(\lambda/\mu)^n$ with $\lambda<\mu, \beta>0$ and $C$ some constant. Using the Euler integral for the hypergeometric function
$$
{}_{2}F_{1}\!\left(\!\begin{array}{c}
1,\beta+n+1 \\[-1.5mm]
n+2
\end{array}; \frac{\lambda}{\mu} \right)=(n+1)\int_0^1\frac{(1-t)^{n}}{\left(1-\frac{\lambda}{\mu}t\right)^{\beta+n+1}}dt,
$$
and the binomial theorem, we can write $T$ as
\begin{align*}
T&=\frac{\beta^2(\lambda/\mu)^2}{\left[\left(1-\lambda/\mu\right)^{-\beta}-1\right]^2}\sum_{n=0}^\infty\frac{(\beta+1)_n(\lambda/\mu)^n}{n!}\int_0^1\int_0^1\frac{(1-t)^n(1-s)^n\,dt\,ds}{\left(1-\frac{\lambda}{\mu}t\right)^{\beta+n+1}\left(1-\frac{\lambda}{\mu}s\right)^{\beta+n+1}}\\
&=\frac{\beta^2(\lambda/\mu)^2}{\left[\left(1-\lambda/\mu\right)^{-\beta}-1\right]^2}\int_0^1\int_0^1\frac{ds\,dt}{\left[\left(1-\frac{\lambda}{\mu}t\right)\left(1-\frac{\lambda}{\mu}s\right)\right]^{\beta+1}}\sum_{n=0}^\infty\frac{(\beta+1)_n}{n!}\left(\frac{\lambda(1-s)(1-t)}{\mu\left(1-\frac{\lambda}{\mu}t\right)\left(1-\frac{\lambda}{\mu}s\right)}\right)^n\\
&=\frac{\beta^2(\lambda/\mu)^2\left(1-\lambda/\mu\right)^{-\beta-1}}{\left[\left(1-\lambda/\mu\right)^{-\beta}-1\right]^2}\int_0^1\int_0^1\frac{ds\,dt}{\left(1-\frac{\lambda}{\mu}ts\right)^{\beta+1}}.
\end{align*}
Using again the binomial theorem for the previous double integral gives a representation of $T$ in terms of the ${}_{3}F_{2}$ hypergeometric function. Therefore
\begin{equation*}
T=\frac{\beta^2(\lambda/\mu)^2\left(1-\lambda/\mu\right)^{-\beta-1}}{\left[\left(1-\lambda/\mu\right)^{-\beta}-1\right]^2}\,{}_{3}F_{2}\!\left(\!\begin{array}{c}
1,1,\beta+1 \\[-1.5mm]
2,2
\end{array}; \frac{\lambda}{\mu} \right).
\end{equation*}
According to \eqref{condBAB}, the absorption rate \(\tilde{\mu}_0\) must satisfy
\begin{equation*}
0\leq\tilde\mu_0 \le\frac{\mu}{\beta}\frac{\left(\dfrac{\mu}{\lambda}-1\right)\left[\left(1-\dfrac{\lambda}{\mu}\right)^{-\beta}-1\right]}{{}_{3}F_{2}\!\left(\!\begin{array}{c}
1,1,\beta+1 \\[-1.5mm]
2,2
\end{array}; \D\frac{\lambda}{\mu} \right)}.
\end{equation*}
Let us now derive further information about the quantities \(T\) and \(\frac{1}{T}\left(\lambda\beta+\frac{1}{A}\right)\) introduced above in terms of the parameter \(z=\lambda/\mu\) and \(\beta>0\). To this end, let us denote in this subsectoin $T=T_\beta(z)$ and $R_\beta(z)=\frac{1}{T}\left(\lambda\beta+\frac{1}{A}\right)$. 

\begin{figure}[htbp]
\centering
\begin{subfigure}[t]{0.49\textwidth}
    \centering
    \includegraphics[width=\textwidth]{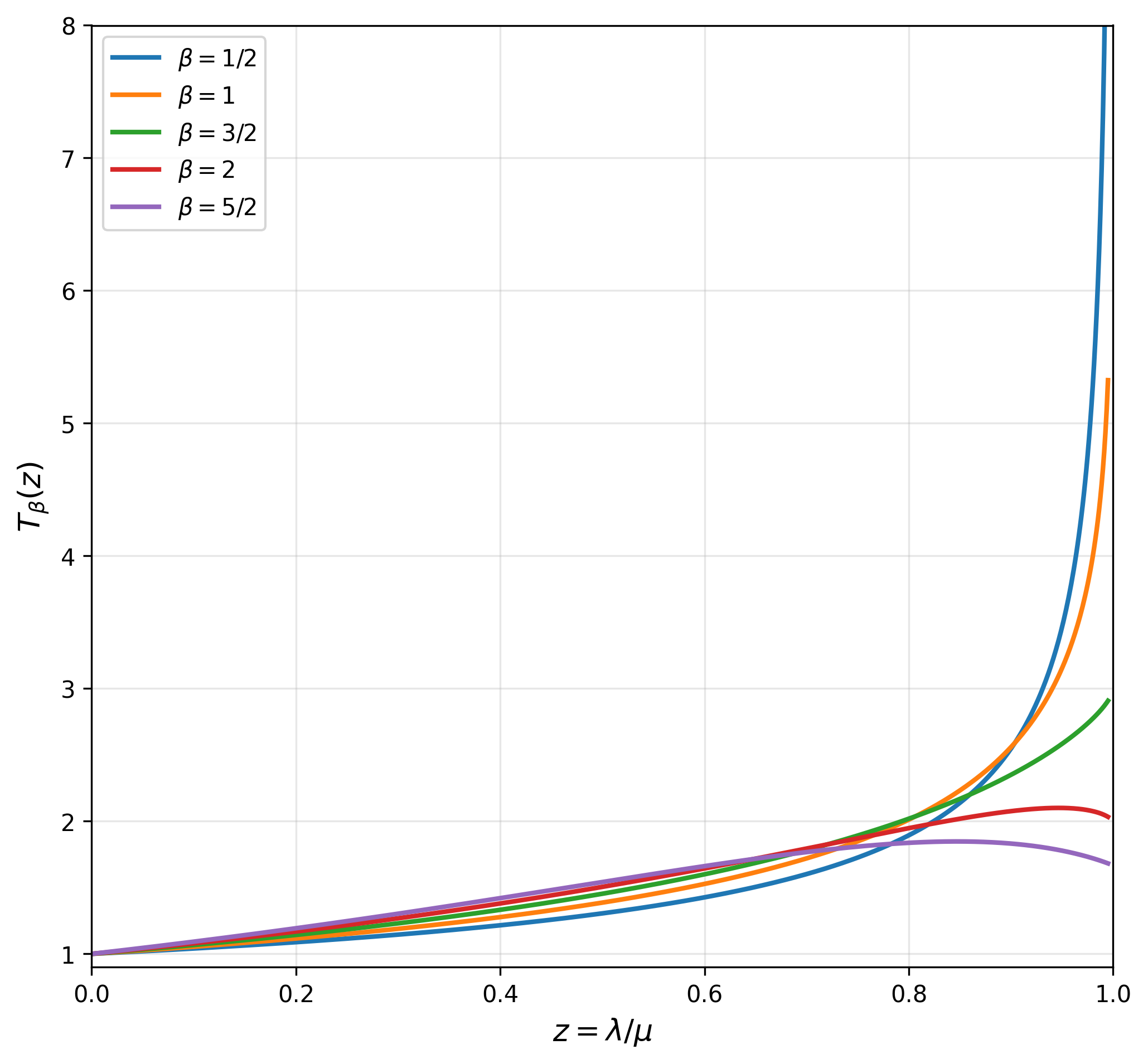}\end{subfigure}
\hfill
\begin{subfigure}[t]{0.49\textwidth}
    \centering
    \includegraphics[width=\textwidth]{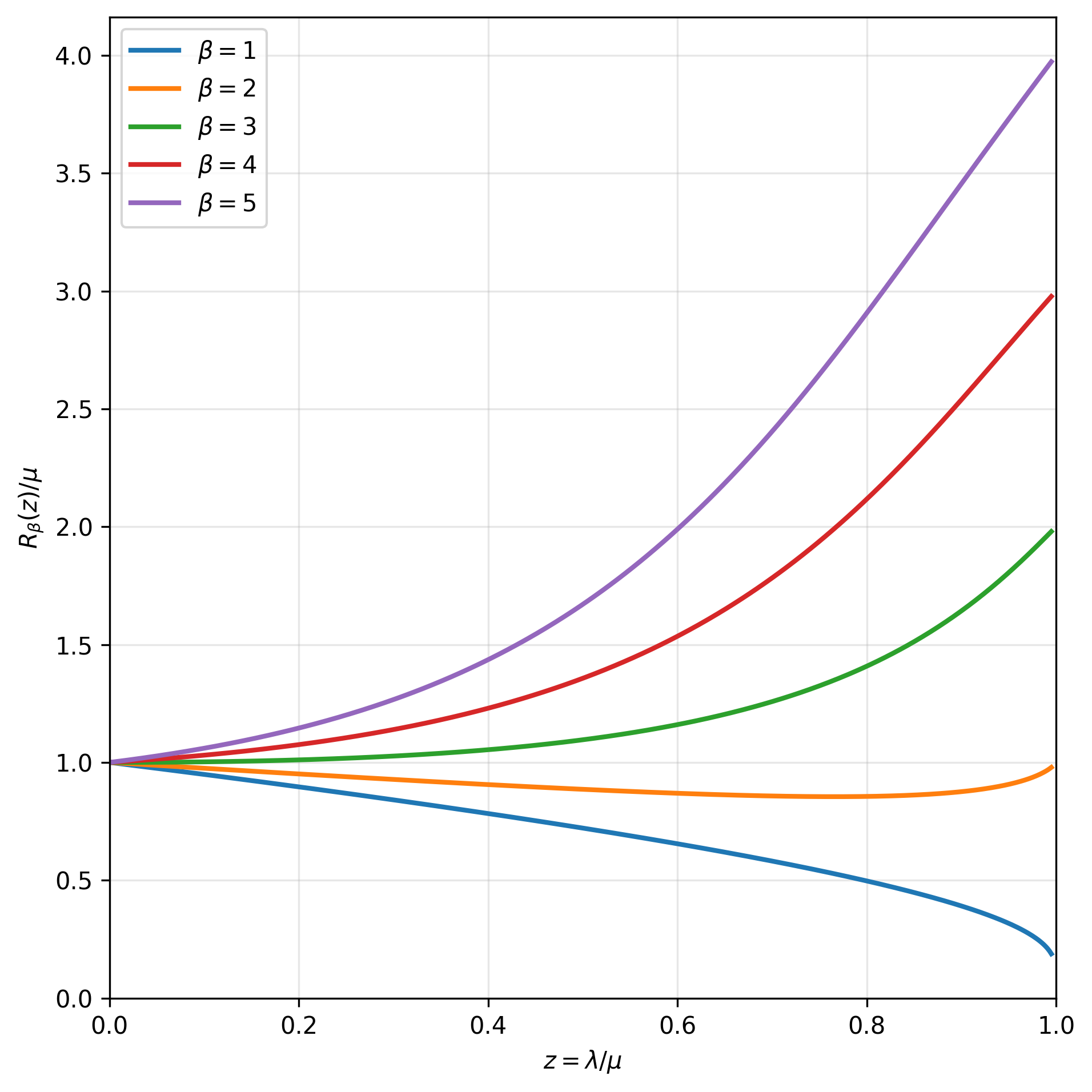}
\end{subfigure}
\caption{Graphs of the functions $T_\beta(z)$ and $R_\beta(z)/\mu$ as functions of $z=\lambda/\mu$ for several values of $\beta$.}
\label{fig:comparison2}
\end{figure}

On one hand, from the Taylor expansion of $T_\beta(z)$ around $z=0$, we obtain
$
T_\beta(z)
=
1+\frac{\beta+1}{4}z
+\mathcal{O}(z^2),
z\to0.
$
Therefore we get $\lim_{z\to0^+}T_\beta(z)=1$ and $\lim_{z\to0^+}T_\beta'(z)=\frac{\beta+1}{4}>0$. That means that $T_\beta$ is always increasing in a right-neighborhood of 0. Additionally, we have that the asymptotic behavior of \(T_\beta(z)\) as \(z\to1^{-}\) is given, for \(0<\beta\le 1\), by
\[
T_\beta(z)\sim
\begin{cases}
C_{\beta}(1-z)^{\beta-1}, & \mbox{if}\quad 0<\beta<1,\\[1em]
\D\frac{\log(1/(1-z))}{z}, & \mbox{if}\quad \beta=1,
\end{cases}
\]
for certain constant \(C_{\beta}\) depending only on \(\beta\). Therefore,
\[
\lim_{z\to1^-}T_\beta(z)
=
\begin{cases}
\infty, & \mbox{if}\quad 0<\beta\leq1,\\[1em]
\D\frac{\beta}{\beta-1}, & \mbox{if}\quad \beta>1.
\end{cases}
\]
In contrast with the case $\lambda>\mu$ (see Section \ref{sec42}), $T_\beta$ is not monotone in general (see Figure~\ref{fig:comparison2}). Finally, for fixed $0<z<1$ we have that $T_\beta(z)\to1$ as $\beta\to\infty$. Indeed, for fixed $z\in(0,1)$, a Laplace-type analysis of the integral representation
$
{}_3F_2(\beta+1,1,1;2,2;z)
=
\int_0^1\int_0^1(1-zst)^{-\beta-1}\,ds\,dt
$
shows that, as $\beta\to\infty$, the integral concentrates near $(s,t)=(1,1)$ and
$
{}_3F_2(\beta+1,1,1;2,2;z)
\sim
\frac{(1-z)^{1-\beta}}{\beta^2z^2}.
$
Since
$
\frac{\beta^2z^2(1-z)^{-\beta-1}}
{\bigl((1-z)^{-\beta}-1\bigr)^2}
\sim
\beta^2z^2(1-z)^{\beta-1},
$
the two asymptotic factors cancel, yielding
$
T_\beta(z)\to1,
$ as $\beta\to\infty.$

On the other hand, the expansion of $R_\beta(z)/\mu$ around $z=0$ is given by
$
R_\beta(z)/\mu
=
1+\frac{\beta-3}{4}z
+\mathcal{O}(z^2), z\to0.
$
Therefore, we obtain $\lim_{z\to0^+}R_\beta(z)/\mu=1$ and $\lim_{z\to0^+}R_\beta'(z)/\mu=\frac{\beta-3}{4}$. Hence \(R_\beta(z)\) is initially decreasing for \(0<\beta<3\), stationary to first order for \(\beta=3\), and initially increasing for \(\beta>3\). Additionally, we always have 
$$
\lim_{z\to1^-}\frac{R_\beta(z)}{\mu}=
\begin{cases}
0, & \mbox{if}\quad 0<\beta\leq1,\\[1em]
\beta-1, & \mbox{if}\quad \beta>1.
\end{cases}
$$
Finally, for fixed $0<z<1$, we have that $R_\beta(z)/\beta\to \lambda$ as $\beta\to\infty$. Indeed, as $\beta\to\infty$, we have
$
\frac1A=
\frac{\beta\lambda(1-z)^\beta}{1-(1-z)^\beta}
=\beta\lambda\mathcal{O}((1-z)^\beta).
$
Hence
$
R_\beta(z)
=\frac{\beta\lambda+1/A}{T_\beta(z)}
\sim\beta\lambda.
$
Therefore, for every fixed $z\in(0,1)$, $R_\beta(z)/\beta\to\lambda=\mu z$ as $\beta\to\infty$. Notice, however, that the asymptotic equivalence $R_\beta(z)\sim\beta\lambda$ is not uniform as $z\to0^+$, since $R_\beta(0)=\mu$ for every $\beta$, whereas $\beta\lambda=\beta\mu z$ vanishes at $z=0$.

\medskip

The corresponding graphs of \(T_\beta(z)\) and \(R_\beta(z)\) for several representative values of \(\beta\) are displayed in Figure~\ref{fig:comparison2}.

\medskip

As for the Darboux spectral measure, the original spectral measure of the dual process is given by \eqref{meixasdf}. We must therefore consider the Geronimus transformation of \eqref{meixasdf}, namely \eqref{spmec2e}. From \eqref{mmm1} it follows that
$$
m_{-1}=\frac{1}{\lambda\beta}\left[1-\left(1-\frac{\lambda}{\mu}\right)^\beta\right].
$$
Substituting this expression into \eqref{spmec2e} and rearranging the resulting series, we obtain after straightforward simplifications
$$
d\widetilde\psi(x)=\frac{x_0+\tilde\mu_0}{\lambda\beta}\left(1-\frac{\lambda}{\mu}\right)^\beta\sum_{n=0}^{\infty}\frac{(\beta)_n}{n!}\left(\frac{\lambda}{\mu}\right)^n\delta_{(\mu-\lambda) n}(dx)+\left(1-\frac{x_0+\tilde\mu_0}{\lambda\beta}\right)\delta_0(dx),
$$
which shows that \(d\widetilde{\psi}\) is an affine combination of the original spectral measure associated with the Meixner polynomials and the Dirac measure at \(x=0\). Equivalently, the masses at the positive support points are rescaled, with the remaining mass concentrated at the origin. As for the orthogonal polynomials, a lengthy but straightforward computation, making use of \cite[(9.10.6)]{KS} and the self-duality of Meixner polynomials together with the property $Q_n(0)=Q_{n-1}(0)+\frac{(\beta)_n}{n!}(\lambda/\mu)^n$ (see \eqref{Qn0mm}), yields $\tilde Q_0=1$ and, for $n\geq1$,
\begin{multline*}
\tilde Q_n(x)=\frac{s_n(\beta)_n}{n!}\left(\frac{\lambda}{\mu}\right)^n\left[M_n\left(\frac{x}{\mu-\lambda};\beta,\frac{\lambda}{\mu}\right)\right.\\ \left.-\left(1-\frac{x_0+\tilde\mu_0}{\lambda\beta}\right)\frac{(\beta+1)_{n-1}(\lambda/\mu)^{n-1}}{(n-1)!u_{n-1}}M_{n-1}\left(\frac{x}{\mu-\lambda}-1;\beta+1,\frac{\lambda}{\mu}\right)\right].
\end{multline*}
It is worth noting that, in the particular case \(x_0+\tilde\mu_0=\lambda\beta\), the above expression reduces to the classical Meixner polynomials.

\subsection{Case $\lambda>\mu$}\label{sec42}

Now we have $\mu/\lambda<1$ and therefore $A<\infty$ and $B=\infty$, that is, the birth--death process is transient. The value of $A$ in \eqref{A_B_linear} is given in terms of the Gauss hypergeometric function ${_2}F_1$. Using Euler's integral representation (see \cite[Theorem 2.2.1]{Ask}) we have that
\begin{equation}\label{AAlin}
A=\frac{1}{\lambda\beta}\,{}_{2}F_{1}\!\left(\!\begin{array}{c}
1,1 \\[-1.5mm]
\beta+1
\end{array}; \frac{\mu}{\lambda} \right)=\frac{1}{\lambda}\int_0^1\frac{(1-t)^{\beta-1}}{1-\frac{\mu}{\lambda}t}dt.
\end{equation}
In general, this expression cannot be simplified further. 
However, when \(\beta = m \in \mathbb{N}\), we obtain the following formula:
\begin{equation}\label{Aforbm}
A=\frac{1}{\lambda-\mu}\left(1-\frac{\lambda}{\mu}\right)^m\left[\log\left(1-\frac{\mu}{\lambda}\right)+\sum_{k=1}^{m-1}\frac{(-1)^k}{k}\binom{m-1}{k}\left(1-\left(1-\frac{\mu}{\lambda}\right)^{-k}\right)\right],
\end{equation}
which can be proved by induction using that, for $m=1$, $A=-\frac{1}{\mu}\log\left(1-\frac{\mu}{\lambda}\right)$ and the following formula
$$
\frac{1}{\lambda\beta}\,{}_{2}F_{1}\!\left(\!\begin{array}{c}
1,1 \\[-1.5mm]
\beta+1
\end{array}; \frac{\mu}{\lambda} \right)=\frac{1-\lambda/\mu}{\lambda(\beta-1)}\,{}_{2}F_{1}\!\left(\!\begin{array}{c}
1,1 \\[-1.5mm]
\beta
\end{array}; \frac{\mu}{\lambda} \right)+\frac{1}{\mu(\beta-1)}.
$$
It is well known that the spectral measure associated with the infinitesimal generator $\mathcal{A}$ is given by
\begin{equation*}
    d\psi(x)=\left(1-\frac{\mu}{\lambda}\right)^{\beta}\sum_{n=0}^{\infty}\frac{(\beta)_n (\mu/\lambda)^n}{n!}\delta_{\xi_n}(dx), \quad \xi_n=(\lambda-\mu)(n+\beta).
\end{equation*}
Again, the masses are the negative binomial or Pascal probabilities, but now with parameters $r=\beta$ and $p=1-\mu/\lambda$. The corresponding orthogonal polynomials can be written in terms of the Meixner polynomials $M_n(x;b,c),n\geq0$. Indeed,
\begin{equation*}
    Q_n(x)=\left(\frac{\mu}{\lambda}\right)^nM_n\left(\frac{x}{\lambda-\mu}-\beta;\beta,\frac{\mu}{\lambda}\right),\quad n\geq0.
\end{equation*}
Now, however, the birth--death process is transient, so there is no stationary distribution. Since $B=\infty$ in this case, we cannot use formula \eqref{posss2}. However, we can obtain an expression for the partial sums in terms of the Gauss hypergeometric function. In particular, we have
\begin{equation}\label{xpsm}
\frac{1}{\pi_n}\sum_{k=0}^n\pi_k={}_{2}F_{1}\!\left(\!\begin{array}{c}
1,-n \\[-1.5mm]
1-\beta-n
\end{array}; \frac{\mu}{\lambda} \right).
\end{equation}
Also, for the partial sums of $A$, since $A<\infty$, we may obtain the following expression
\begin{equation}\label{xasm}
\sum_{k=0}^{n-1}\frac{1}{\lambda_k\pi_k}
=\sum_{k=0}^{n-1}\frac{k!(\mu/\lambda)^k}{\lambda(k+\beta)(\beta)_k}=\frac{1}{\lambda\beta}\,{}_{2}F_{1}\!\left(\!\begin{array}{c}
1,1 \\[-1.5mm]
\beta+1
\end{array}; \frac{\mu}{\lambda} \right)-\frac{1}{\lambda_n\pi_n}\,{}_{2}F_{1}\!\left(\!\begin{array}{c}
1,n+1 \\[-1.5mm]
\beta+n+1
\end{array}; \frac{\mu}{\lambda} \right).
\end{equation}

\medskip

\hspace{-.35cm}\underline{\textbf{LU factorization}}. The coefficients $\tilde s_n, \tilde r_n, \tilde x_n,$ and $\tilde y_n$ in \eqref{ssnn} and \eqref{rrnn}--\eqref{yynn} can be explicitly computed. Given that \(\mu_0=0\), we have that $\tilde q_n=-1$ and $\tilde t_n=1-\frac{\hat\mu_0}{\lambda\beta}, \tilde t_0=1$. Therefore, using \eqref{xpsm}, we have
\begin{equation}\label{srxymmi3}
\begin{split}
\tilde s_n&=\frac{\hat\mu_0\mu^nn!}{\beta(\beta)_n\lambda^{n+1}}+\left(1-\frac{\hat\mu_0}{\lambda\beta}\right)\,{}_{2}F_{1}\!\left(\!\begin{array}{c}
1,-n \\[-1.5mm]
1-\beta-n
\end{array}; \frac{\mu}{\lambda} \right),\quad n\geq0,\\
\tilde r_n&=-\frac{\hat\mu_0\mu^nn!}{\beta(\beta)_n\lambda^{n+1}}-\left(1-\frac{\hat\mu_0}{\lambda\beta}\right)\frac{\mu n}{\lambda(n+\beta-1)}\,{}_{2}F_{1}\!\left(\!\begin{array}{c}
1,-n+1 \\[-1.5mm]
2-\beta-n
\end{array}; \frac{\mu}{\lambda} \right),\quad n\geq1,\\
\tilde x_n&=\frac{\lambda(n+\beta)}{\tilde s_n}=-\tilde y_n,\quad n\geq0.
\end{split}
\end{equation}
With this, we can compute the coefficients $\hat\lambda_n,\hat\mu_{n+1},n\geq0,$ in \eqref{newbdcoe}, given by
$$
\hat\lambda_n=\lambda(n+\beta)\frac{\tilde s_{n+1}}{\tilde s_n},\quad\hat\mu_{n+1}=\frac{\mu(n+\beta+1)(n+1)}{n+\beta}\frac{\tilde s_n}{\tilde s_{n+1}},\quad n\geq0.
$$
In order to obtain the Darboux transformation \(\widehat{\mathcal{A}}\) as the infinitesimal generator of a new birth--death process, we apply the second part of Theorem \ref{theom33}, since $\mu_0=0$. Moreover, since $B=\infty$, the parameter \(\hat\mu_0\) must satisfy condition~\eqref{condLU2}, namely,
\[
0\le\hat\mu_0\le\lambda\beta.
\]
As for the Darboux spectral measure, since \(\tilde y_0=-\lambda\beta\), we have
\begin{align*}
d\widehat\psi(x)&=\frac{x}{\lambda\beta}d\psi(x)=\left(1-\frac{\mu}{\lambda}\right)^{\beta}\sum_{n=0}^\infty\frac{(\lambda-\mu)(n+\beta)(\beta)_n(\mu/\lambda)^n}{\lambda\beta n!}\delta_{(\lambda-\mu)(n+\beta)}(dx)\\
&=\left(1-\frac{\mu}{\lambda}\right)^{\beta+1}\sum_{n=0}^{\infty}\frac{(\beta+1)_n (\mu/\lambda)^n}{n!}\delta_{(\lambda-\mu)(n+\beta)}(dx),
\end{align*}
that is, the transformed measure has the same support as $d\psi$, while its weights correspond to those of a negative binomial distribution with parameter $\beta+1$. Similarly, the Darboux birth--death polynomials \((\hat Q_n)_{n\in\mathbb{N}_0}\) are given by
\begin{equation*}
        \hat Q_n(x)=\frac{\tilde x_n}{\lambda\beta}\left(\frac{\mu}{\lambda}\right)^nM_n\left(\frac{x}{\lambda-\mu}-\beta;\beta+1,\frac{\mu}{\lambda}\right),\quad n\geq0.
\end{equation*}
Substituting $\hat\mu_0=0$ in \eqref{srxymmi3} gives directly the coefficients $\tilde s_n, \tilde r_n, \tilde x_n,$ and $\tilde y_n$. For the other limiting case, that is, $\hat\mu_0=\lambda\beta$, we get the dual process, see Remark \ref{remdual}.

\medskip

\hspace{-.35cm}\underline{\textbf{UL factorization}}. The coefficients $x_n, y_n, s_n,$ and $r_n$ in \eqref{xxnn2} and \eqref{yynn2}--\eqref{rrnn2} can be explicitly computed from the sequence $(u_n)_{n \in \mathbb{N}_0}$, which in this case is given, using \eqref{xasm}, by
\begin{equation}\label{uunlin}
u_n=1-(x_0+\tilde\mu_0)\left[\frac{1}{\lambda\beta}\,{}_{2}F_{1}\!\left(\!\begin{array}{c}
1,1 \\[-1.5mm]
\beta+1
\end{array}; \frac{\mu}{\lambda} \right)-\frac{1}{\lambda_n\pi_n}\,{}_{2}F_{1}\!\left(\!\begin{array}{c}
1,n+1 \\[-1.5mm]
\beta+n+1
\end{array}; \frac{\mu}{\lambda} \right)\right],\quad n\geq0.
\end{equation}
The probabilistic interpretation of these coefficients was discussed in \cite[Section~4.1]{AdI}. When $\tilde\mu_0=0$, Remark~\ref{remsimpl} gives the rates $\tilde\lambda_0=x_0$, $\tilde\lambda_n=\mu_n u_{n-1}/u_n$ and $\tilde\mu_{n+1}=\lambda_n u_{n+1}/u_n$; for $\tilde\mu_0>0$ they follow from \eqref{newbdcoe2}--\eqref{laconcha3r}.

Recall that, in this context, we have to apply Theorem \ref{theom37}. Since $A<\infty$ and $B=\infty$, if $x_0+\tilde\mu_0<1/A$, then $T=\infty$, which implies that $\tilde\mu_0=0$. Therefore the free parameter \(x_0\) must satisfy (see \eqref{cond1UL})
\begin{equation*}
0 < x_0 < \frac{\lambda\beta}{{}_{2}F_{1}\!\left(\!\begin{array}{c}
1,1 \\[-1.5mm]
\beta+1
\end{array}; \dfrac{\mu}{\lambda} \right)
}.
\end{equation*}
On the other hand, if $x_0+\tilde\mu_0=1/A$, then the series $T$ can be written as
$$
T=\frac{\beta^2}{{}_{2}F_{1}\!\left(\!\begin{array}{c}
1,1 \\[-1.5mm]
\beta+1
\end{array}; \dfrac{\mu}{\lambda} \right)^2}\sum_{n=0}^\infty\frac{n!}{(n+\beta)(\beta)_{n+1}}\left(\frac{\mu}{\lambda}\right)^n\,{}_{2}F_{1}\!\left(\!\begin{array}{c}
1,n+1 \\[-1.5mm]
\beta+n+1
\end{array}; \frac{\mu}{\lambda} \right)^2.
$$
Using the Euler integral for the hypergeometric function
$$
{}_{2}F_{1}\!\left(\!\begin{array}{c}
1,n+1 \\[-1.5mm]
\beta+n+1
\end{array}; \frac{\mu}{\lambda} \right)=\frac{(\beta)_{n+1}}{n!}\int_0^1\frac{t^n(1-t)^{\beta-1}}{1-\frac{\mu}{\lambda}t}dt=\frac{(\beta)_{n+1}}{n!}I_n,
$$
we can write $T$, using \eqref{AAlin}, as
$$
T=\frac{1}{\lambda^2A^2}\sum_{n=0}^\infty\left(\frac{\mu}{\lambda}\right)^n\frac{(\beta)_{n}}{n!}I_n^2.
$$
The series $T$ actually converges, since $\mu<\lambda$, and $\frac{(\beta)_n}{n!}\sim\frac{n^{\beta-1}}{\Gamma(\beta)}$, $I_n\sim\frac{\Gamma(\beta)}{1-\frac{\mu}{\lambda}}n^{-\beta}$, as $n\to\infty$. Using the integral representations, $T$ can be written as 
$$
T=\frac{\D\iint_{[0,1]^2}\dfrac{(1-t)^{\beta-1}(1-s)^{\beta-1}}{\left(1-\frac{\mu}{\lambda}t\right)\left(1-\frac{\mu}{\lambda}s\right)\left(1-\frac{\mu}{\lambda}st\right)^\beta}dtds}{\left(\D\int_0^1\dfrac{(1-t)^{\beta-1}}{1-\frac{\mu}{\lambda}t}dt\right)^2}.
$$
Therefore, following \eqref{condBAB}, the absorption rate $\tilde\mu_0$ can be chosen in the following range:
$$
0\leq\tilde\mu_0\leq\frac{1}{AT},
$$
which can be written, using the integral representations, as
$$
0\leq\tilde\mu_0\leq\lambda\,\frac{\D\int_0^1\dfrac{(1-t)^{\beta-1}}{1-\frac{\mu}{\lambda}t}dt}{\D\iint_{[0,1]^2}\dfrac{(1-t)^{\beta-1}(1-s)^{\beta-1}}{\left(1-\frac{\mu}{\lambda}t\right)\left(1-\frac{\mu}{\lambda}s\right)\left(1-\frac{\mu}{\lambda}st\right)^\beta}dtds}.
$$
Let us now derive further information about the quantities \(T\) and \((AT)^{-1}\) introduced above in terms of the parameter \(z=\mu/\lambda\) and \(\beta>0\). To this end, let \(N_\beta(z)\) and \(D_\beta^2(z)\) denote the numerator and denominator of \(T\), respectively. Then
\[
T_\beta(z):=T=\frac{N_\beta(z)}{D_\beta^2(z)},\quad
R_\beta(z):=\lambda\frac{D_\beta(z)}{N_\beta(z)},\quad
z=\frac{\mu}{\lambda}.
\]
Since \(0<\mu/\lambda<1\), both functions are defined on the interval \((0,1)\). It is straightforward to verify that
\(
\lim_{z\to0^+}D_\beta(z)=\frac{1}{\beta},
\lim_{z\to0^+}N_\beta(z)=\frac{1}{\beta^2},
\)
and therefore
\(
\lim_{z\to0^+}T_\beta(z)=1.
\)
Expanding \((1-zt)^{-1}\) around \(z=0\), one obtains
\(
T_\beta(z)
=
1+\frac{\beta}{(\beta+1)^2}z+\mathcal{O}(z^2), z\to0^+.
\)
In particular, \(T_\beta'(0)>0\), showing that \(T_\beta(z)\) starts at \(1\) and is initially increasing. In fact, it is possible to prove that $T_\beta(z)$ is strictly increasing on $(0,1)$. This follows writing \(
w_z(t)=\frac{(1-t)^{\beta-1}}{1-zt}\) and \(
K_z(t,s)=(1-zts)^{-\beta},
\)
so that
\(
D_\beta(z)=\int_0^1 w_z(t)\,dt\) and \(
N_\beta(z)=\int_0^1\int_0^1
w_z(t)w_z(s)K_z(t,s)\,dt\,ds.
\)
Since
\(
\partial_z w_z(t)=\frac{t}{1-zt}w_z(t)\) and \(
\partial_z K_z(t,s)
=
\frac{\beta ts}{1-zts}K_z(t,s),
\)
we obtain, by symmetry,
\[
N_\beta'(z)=
2\,
\int_0^1\int_0^1
\frac{t}{1-zt}
w_z(t)w_z(s)K_z(t,s)\,dt\,ds
+
\int_0^1\int_0^1
\frac{\beta ts}{1-zts}
w_z(t)w_z(s)K_z(t,s)\,dt\,ds.
\]
For each fixed $s>0$, both
\(
\frac{t}{1-zt}\) and \(
K_z(t,s)
\)
are strictly increasing in $t$. Hence the weighted Chebyshev
integral inequality (see \cite[Chapter IX]{MPF} for $f(t)=\frac{t}{1-zt}, g(t)=K_z(t,s)$ and $p(t)=w_z(t)$) gives
\[
\frac{
\displaystyle
\int_0^1
\frac{t}{1-zt}K_z(t,s)w_z(t)\,dt
}{
\displaystyle
\int_0^1K_z(t,s)w_z(t)\,dt
}
>
\frac{
\displaystyle
\int_0^1
\frac{t}{1-zt}w_z(t)\,dt
}{
D_\beta(z)
}
=
\frac{D_\beta'(z)}{D_\beta(z)}.
\]
After integrating with respect to $s$, it follows that
\(
\frac{N_\beta'(z)}{N_\beta(z)}
>
2\frac{D_\beta'(z)}{D_\beta(z)}.
\)
Therefore
\(
\frac{T_\beta'(z)}{T_\beta(z)}
=
\frac{N_\beta'(z)}{N_\beta(z)}
-
2\frac{D_\beta'(z)}{D_\beta(z)}
>0,
\)
and hence $T_\beta(z)$ is strictly increasing on $(0,1)$. The only singularity on $[0,1]$ is at the endpoint $z=1$. The asymptotic behavior of \(T_\beta(z)\) as \(z\to1^{-}\) is given by
\[
T_\beta(z)\sim
\begin{cases}
C_{\beta,1}(1-z)^{-\beta}, & \mbox{if}\quad 0<\beta<1,\\[1em]
\D\frac{C_{\beta,2}}{(1-z)\log^2(1-z)}, & \mbox{if}\quad \beta=1,\\[1em]
C_{\beta,3}(1-z)^{\beta-2}, & \mbox{if}\quad 1<\beta<2,\\[1em]
-\log(1-z), & \mbox{if}\quad \beta=2,
\end{cases}
\]
for certain constants \(C_{\beta,i}\), \(i=1,2,3\), depending only on \(\beta\). Hence,
\[
\lim_{z\to1^-}T_\beta(z)
=
\begin{cases}
\infty, & \mbox{if}\quad 0<\beta\leq2,\\[1em]
\D\frac{\beta-1}{\beta-2}, & \mbox{if}\quad \beta>2.
\end{cases}
\]
Finally, for fixed $0<z<1$ we have that $T_\beta(z)\to1$ as $\beta\to\infty$. Indeed, for fixed $z\in(0,1)$, as $\beta\to\infty$ the factors
$(1-t)^{\beta-1}$ and $(1-s)^{\beta-1}$ concentrate near $t=s=0$. Expanding the remaining factors around this point yields
\(
D_\beta(z)
=\frac{1}{\beta}\left(1+\frac{z}{\beta}
+\mathcal{O}(\beta^{-2})\right),
N_\beta(z)
=\frac{1}{\beta^2}\left(1+\frac{3z}{\beta}
+\mathcal{O}(\beta^{-2})\right).
\)
Consequently,
\(
T_\beta(z)
=\frac{N_\beta(z)}{D_\beta^2(z)}
=1+\frac{z}{\beta}+\mathcal{O}(\beta^{-2}),
\)
and hence $T_\beta(z)\to1$ as $\beta\to\infty$.

A similar analysis can be carried out for \(R_\beta(z)\). In this case, the expansion of \(R_\beta(z)/\lambda\) around $z=0$ is given by \(R_\beta(z)/\lambda=\beta-\frac{\beta(2\beta+1)}{(\beta+1)^2}z+\mathcal{O}(z^2)\). Therefore
\(
R_\beta(0)=\lambda\beta,
\)
and \(R_\beta(z)\) is initially decreasing since $R_\beta'(0)<0$. 
In fact, for fixed $\lambda$, \(R_\beta(z)\) is strictly decreasing on $(0,1)$ since \(
R_\beta(z)/\lambda
=
\frac{1}{D_\beta(z)T_\beta(z)}.
\)
Since both $D_\beta(z)$ and $T_\beta(z)$ are strictly increasing, $R_\beta(z)/\lambda$ is strictly decreasing on $(0,1)$. Furthermore,
\[
\lim_{z\to1^-}R_\beta(z)
=
\begin{cases}
0, & \mbox{if}\quad 0<\beta\leq2,\\[1em]
\lambda(\beta-2), & \mbox{if}\quad \beta>2.
\end{cases}
\]
Finally, for fixed $0<z<1$ we have that $R_\beta(z)/\beta\to\lambda$ as $\beta\to\infty$. Indeed, using the previous expansions for $D_\beta(z)$ and $N_\beta(z)$ in terms of $\beta$, we obtain \(
R_\beta(z)
=\lambda\frac{D_\beta(z)}{N_\beta(z)}
=\lambda\beta-2\lambda z+\mathcal{O}(\beta^{-1}),
\)
so that $R_\beta(z)\sim\lambda\beta$ and, in particular,
$R_\beta(z)/\beta\to\lambda$ as $\beta\to\infty$.

\begin{figure}[htbp]
\centering
\begin{subfigure}[t]{0.48\textwidth}
    \centering
    \includegraphics[width=\textwidth]{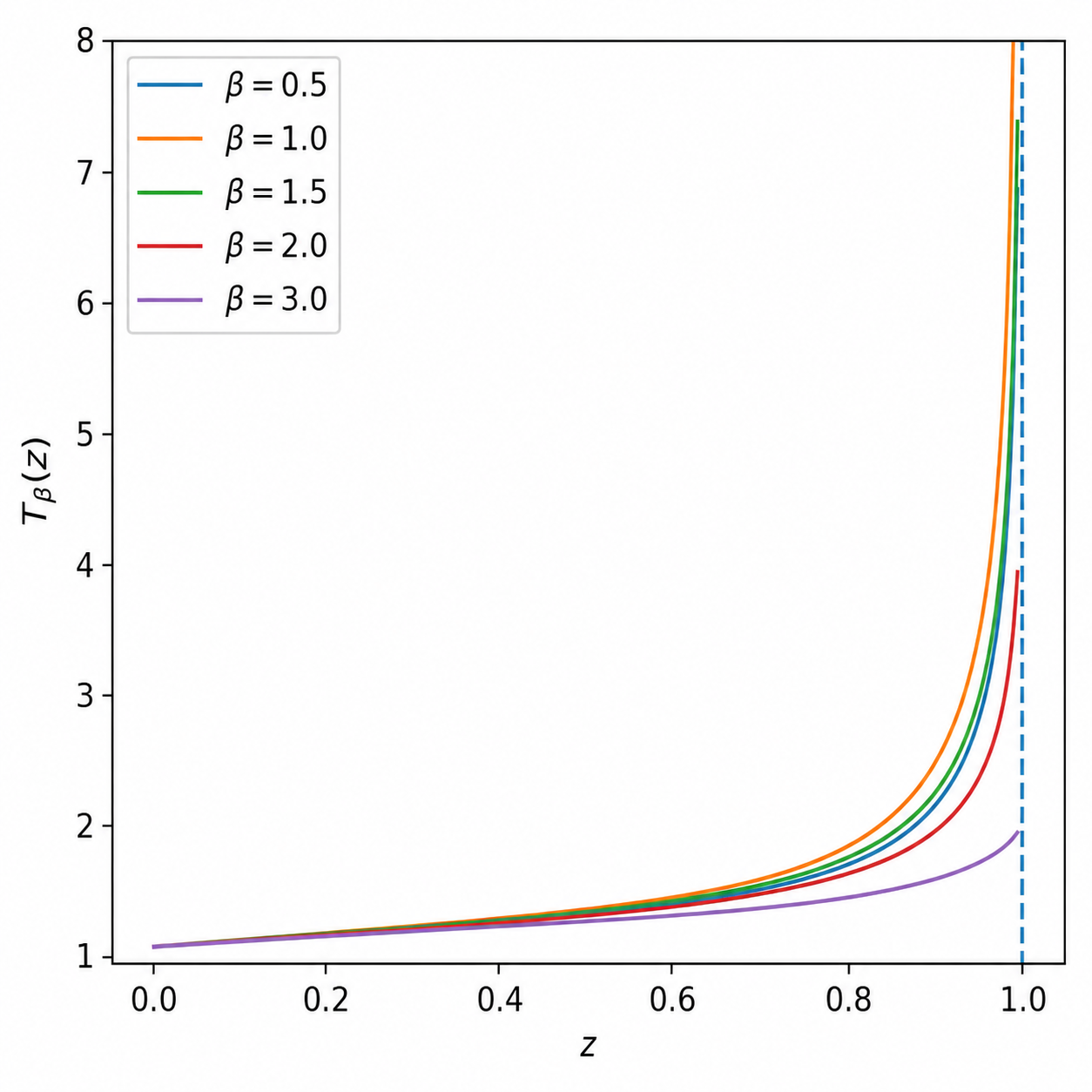}\end{subfigure}
\hfill
\begin{subfigure}[t]{0.48\textwidth}
    \centering
    \includegraphics[width=\textwidth]{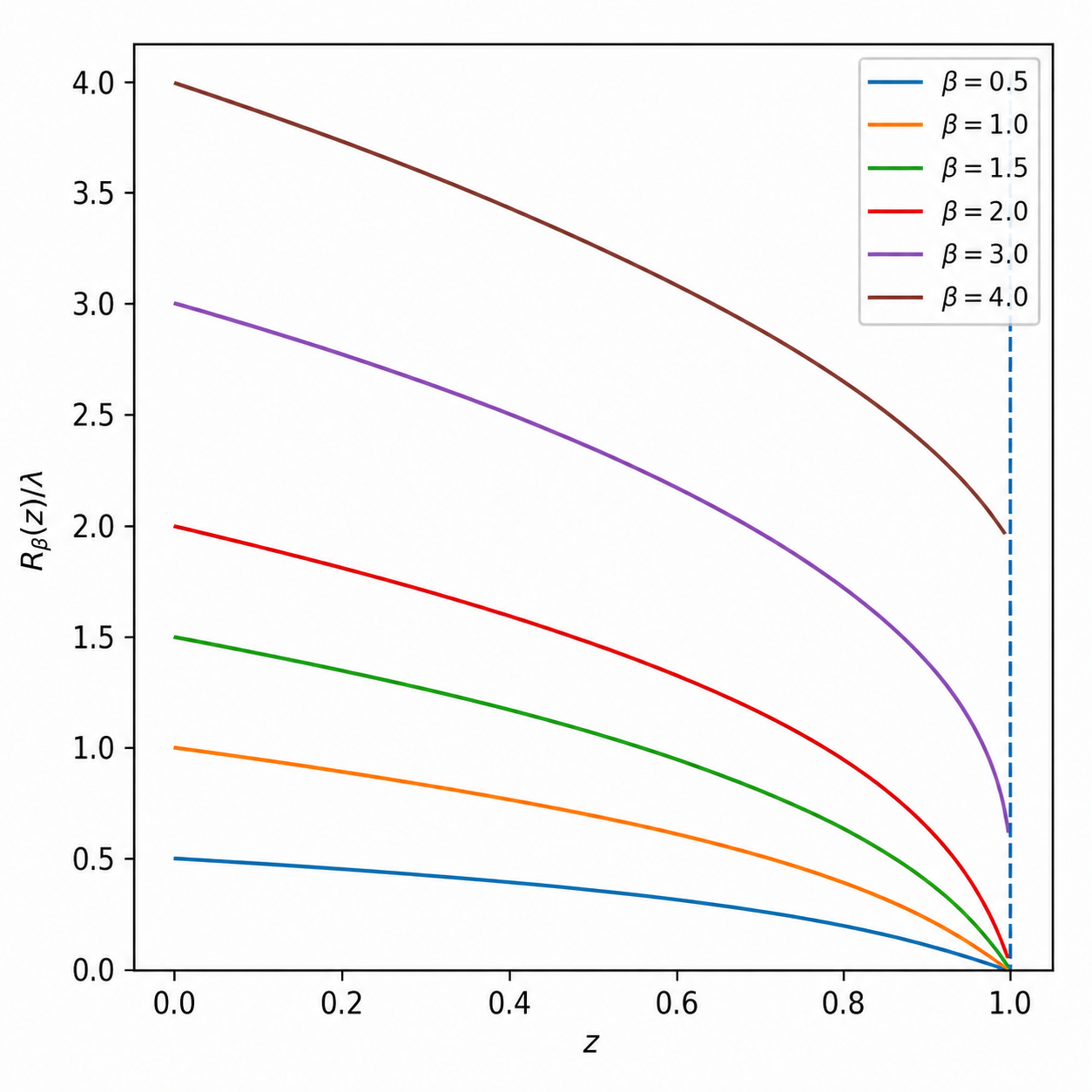}
\end{subfigure}
\caption{Graphs of the functions $T_\beta(z)$ and $R_\beta(z)/\lambda$ as functions of $z=\mu/\lambda$ for several values of $\beta$.}
\label{fig:comparison}
\end{figure}

\medskip
The corresponding graphs of \(T_\beta(z)\) and \(R_\beta(z)\) for several representative values of \(\beta\) are displayed in Figure~\ref{fig:comparison}.

\medskip

The explicit evaluation of the series \(T\) does not appear to be straightforward. However, when \(\beta = m \in \mathbb{N}\), the sum can be expressed in terms of \emph{polylogarithms}. For illustration, we present the case \(\beta = 1\). Using \eqref{Aforbm} for $\beta=1$ and after long but straightforward computations, we obtain
$$
\frac{1}{AT}=\frac{(\lambda-\mu)\log\left(\dfrac{\lambda}{\lambda-\mu}\right)}{\mbox{Li}_2\left(\dfrac{\mu}{\lambda}\right)},
$$
where $\mbox{Li}_2(x)=-\int_0^x\frac{\log(1-u)}{u}du=\sum_{n=1}^\infty\frac{x^n}{n^2}$ is the \emph{dilogarithm}.

\smallskip

Under the conditions considered above, the Darboux transformation $\widetilde{\mathcal{A}}$ becomes the infinitesimal generator of a new birth--death process with birth--death rates $\{\tilde{\lambda}_n, \tilde{\mu}_n\}$ given by \eqref{newbdcoe2} and \eqref{laconcha3r}. The corresponding spectral measure of $\widetilde{\mathcal{A}}$ is described by the Geronimus transform of $d\psi$ in \eqref{spmec2e}. For that, we compute the moment $m_{-1}$ of the spectral measure $d\psi$ using formula \eqref{mmm1}, yielding $m_{-1}=A$. Hence, the transformed spectral measure of the Darboux transformation takes the form
\begin{equation*}
d\widetilde\psi(x)=(x_0+\tilde\mu_0)\left(1-\frac{\mu}{\lambda}\right)^{\beta}\sum_{n=0}^{\infty}\frac{(\beta)_n (\mu/\lambda)^n}{n!(\lambda-\mu)(n+\beta)}\delta_{(\lambda-\mu)(n+\beta)}(dx)+(1-(x_0+\tilde\mu_0)A)\delta_0(dx).
\end{equation*}
The corresponding orthogonal polynomials can be obtained from \eqref{QDTUL} (see also \eqref{QQbb1}); using the following relation 
$$
{}_{2}F_{1}\!\left(\!\begin{array}{c}
1,n \\[-1.5mm]
\beta+n
\end{array}; \frac{\mu}{\lambda} \right)=1+\frac{\mu n}{\lambda(n+\beta)}{}_{2}F_{1}\!\left(\!\begin{array}{c}
1,n+1 \\[-1.5mm]
\beta+n+1
\end{array}; \frac{\mu}{\lambda} \right),
$$
in $u_n$ \eqref{uunlin} and after some straightforward computations, we obtain $\tilde Q_0(x)=1$ and, for $n\geq1$,
$$
\tilde Q_n(x)=s_n\left(\frac{\mu}{\lambda}\right)^{n-1}\left[\frac{x_0+\tilde\mu_0}{u_{n-1}\lambda_{n-1}\pi_{n-1}}M_{n-1}\left(\frac{x}{\lambda-\mu}-\beta;\beta,\frac{\mu}{\lambda}\right)-\frac{x}{\lambda\beta}M_{n-1}\left(\frac{x}{\lambda-\mu}-\beta;\beta+1,\frac{\mu}{\lambda}\right)\right].
$$

\subsection{Case $\lambda=\mu$}
Now, the potential coefficients only depend on $\beta$, that is, 
$$
\pi_n=\frac{(\beta)_n}{n!},\quad n\geq0.
$$
Since $\frac{(\beta)_n}{n!}\sim\frac{n^{\beta-1}}{\Gamma(\beta)}$ as $n\to\infty$, the quantities $A$ and $B$ in \eqref{A_B_linear} are given by
\begin{equation*}
A=\begin{cases}
\infty,&\mbox{if}\quad 0<\beta\leq1\\[.5em]
\dfrac{1}{\lambda(\beta-1)},&\mbox{if}\quad \beta>1\\
\end{cases},\quad B=\infty.
\end{equation*}
Therefore, if $0<\beta\leq1$, the birth--death process is null recurrent, while if $\beta>1$ it is transient. The partial sums of $A$ and $B$ can also be computed explicitly. For $n\geq1$,
\begin{equation}\label{Annnn}
\sum_{k=0}^{n-1}\frac{1}{\lambda_k\pi_k}
=
\begin{cases}
\displaystyle
\frac{1}{\lambda(\beta-1)}
\left[1-\frac{n!}{(\beta)_n}\right],
& \beta\neq1,\\[1.2em]
\displaystyle
\frac{H_n}{\lambda},
& \beta=1,
\end{cases}
\end{equation}
where $H_n=\sum_{k=1}^n k^{-1}$ denotes the $n$-th harmonic number, and 
\begin{equation}\label{Bnnnn}
\sum_{k=0}^{n}\pi_k=\frac{(\beta+1)_n}{n!}\quad\Longrightarrow\quad \frac{1}{\pi_n}\sum_{k=0}^{n}\pi_k=1+\frac{n}{\beta},\quad n\geq0.
\end{equation}
It is well known that the spectral measure associated with the infinitesimal generator $\mathcal{A}$ is given by the density function of the \emph{Gamma distribution} with shape parameter $\beta$ and rate parameter $1/\lambda$ (or scale parameter $\lambda$), that is,
\begin{equation*}
    d\psi(x)=\frac{1}{\Gamma(\beta)\lambda^{\beta}}x^{\beta-1}e^{-x/\lambda}dx,\quad x>0.
\end{equation*}
The corresponding orthogonal polynomials can be written in terms of the \emph{Laguerre polynomials} $L_n^{(\alpha)}(x),n\geq0$. Indeed,
\begin{equation*}
    Q_n(x)=\frac{n!}{(\beta)_n}L_{n}^{(\beta-1)}\left(\frac{x}{\lambda}\right),\quad n\geq0.
\end{equation*}

\medskip

\hspace{-.35cm}\underline{\textbf{LU factorization}}. The coefficients $\tilde s_n, \tilde r_n, \tilde x_n,$ and $\tilde y_n$ in \eqref{ssnn} and \eqref{rrnn}--\eqref{yynn} can be explicitly computed. Given that \(\mu_0=0\), we have that $\tilde q_n=-1$ and $\tilde t_n=1-\frac{\hat\mu_0}{\lambda\beta}, \tilde t_0=1$. Therefore, using \eqref{Bnnnn}, we have
\begin{equation}\label{srxymmi4}
\begin{split}
\tilde s_n&=\frac{\hat\mu_0n!}{\lambda\beta(\beta)_n}+\left(1-\frac{\hat\mu_0}{\lambda\beta}\right)\left(1+\frac{n}{\beta}\right),\quad n\geq0,\\
\tilde r_n&=-\frac{\hat\mu_0n!}{\lambda\beta(\beta)_n}-\left(1-\frac{\hat\mu_0}{\lambda\beta}\right)\frac{n}{\beta},\quad n\geq1,\\
\tilde x_n&=\frac{\lambda(n+\beta)}{\tilde s_n}=-\tilde y_n,\quad n\geq0.
\end{split}
\end{equation}
With this, we can compute the coefficients $\hat\lambda_n,\hat\mu_{n+1},n\geq0,$ in \eqref{newbdcoe}, given by
$$
\hat\lambda_n=\lambda(n+\beta)\frac{\tilde s_{n+1}}{\tilde s_n},\quad\hat\mu_{n+1}=\frac{\lambda(n+\beta+1)(n+1)}{n+\beta}\frac{\tilde s_n}{\tilde s_{n+1}},\quad n\geq0.
$$
In order to obtain the Darboux transformation \(\widehat{\mathcal{A}}\) as the infinitesimal generator of a new birth--death process, we apply the second part of Theorem \ref{theom33}, since $\mu_0=0$. Moreover, since $B=\infty$, the parameter \(\hat\mu_0\) must satisfy condition~\eqref{condLU2}, namely,
\[
0\le\hat\mu_0\le\lambda\beta.
\]
As for the Darboux spectral measure, since \(\tilde y_0=-\lambda\beta\), we have
$$
d\widehat\psi(x)=\frac{x}{\lambda\beta}d\psi(x)=\frac{1}{\Gamma(\beta+1)\lambda^{\beta+1}}x^{\beta}e^{-x/\lambda}dx,\quad x>0,
$$
that is, the same Gamma measure \(d\psi\) with the shape parameter shifted from $\beta$ to $\beta+1$. Similarly, the Darboux birth--death polynomials \((\hat Q_n)_{n\in\mathbb{N}_0}\) are given by
    \begin{equation*}
        \hat Q_n(x)=\frac{\tilde x_n n!}{\lambda(\beta)_{n+1}}L_n^{(\beta)}\left(\frac{x}{\lambda}\right),\quad n\geq0.
\end{equation*}
Substituting $\hat\mu_0=0$ in \eqref{srxymmi4} gives
 \begin{equation*}
        \tilde s_n=1+\frac{n}{\beta},  \quad \tilde r_n=-\frac{n}{\beta},\quad \tilde x_n=\lambda\beta=-\tilde y_n,
    \end{equation*}
    where $\tilde s_n,\tilde x_n,\tilde y_n$ are defined for $n\geq0$, whereas $\tilde r_n$ is defined for $n\geq1$. The birth--death rates for the Darboux-transformed process are given by 
    \begin{equation*}
    \hat\lambda_n=\lambda(n+\beta+1),\quad\hat\mu_{n+1}=\lambda(n+1),\quad n\ge0.
\end{equation*}
Note that the Darboux transform is again a linear birth--death process, but with \(\beta\) replaced by \(\beta+1\). For the other limiting case, that is, $\hat\mu_0=\lambda\beta$, we get the dual process, and the birth--death rates for the Darboux birth--death process are given by $\hat{\lambda}_n=\lambda(n+1)$ and $\hat{\mu}_{n}=\lambda(n+\beta)$ for $n\geq0$.

\medskip

\hspace{-.35cm}\underline{\textbf{UL factorization}}. Since $\mu_0=0$, a direct UL factorization of $\mathcal A$ requires $\beta>1$. Otherwise \(A=\infty\), and admissibility would force \(x_0+\tilde\mu_0\le 0\), which is not possible. In this situation we have that $A=\frac{1}{\lambda(\beta-1)}$, that is, the birth--death process is transient. The coefficients $x_n, y_n, s_n,$ and $r_n$ in \eqref{xxnn2} and \eqref{yynn2}--\eqref{rrnn2} can be explicitly computed from the sequence $(u_n)_{n \in \mathbb{N}_0}$, which in this case is given, using \eqref{Annnn}, by
$$
u_n=1-\frac{x_0+\tilde\mu_0}{\lambda(\beta-1)}\left[1-\frac{n!}{(\beta)_n}\right],\quad n\geq0.
$$
The probabilistic interpretation of these coefficients was discussed in \cite[Section~4.1]{AdI}. In general, the resulting expressions are cumbersome, except for certain particular cases.

Recall that, in this context, we have to apply Theorem \ref{theom37}. Since $A<\infty$ and $B=\infty$, if $x_0+\tilde\mu_0<\lambda(\beta-1)$, then $T=\infty$, which implies that $\tilde\mu_0=0$. Therefore the free parameter \(x_0\) must satisfy (see \eqref{cond1UL})
\begin{equation*}
0 < x_0  < \lambda(\beta-1).
\end{equation*}
On the other hand, if $x_0+\tilde\mu_0=\lambda(\beta-1)$, then $u_n=\pi_n^{-1},n\geq0.$ Since $\frac{(\beta)_n}{n!}\sim\frac{n^{\beta-1}}{\Gamma(\beta)}$ as $n\to\infty$, we have that the series $T$ is given by
$$
T=\sum_{n=0}^\infty \frac{1}{\pi_n}=\sum_{n=0}^\infty \frac{n!}{(\beta)_n}=\begin{cases}\infty,&\mbox{if}\quad 1<\beta\leq2,\\[.5em]
\dfrac{\beta-1}{\beta-2},&\mbox{if}\quad \beta>2.\end{cases}
$$
Therefore we have two situations:
\begin{itemize}[leftmargin=0.25in]
\item If $1<\beta\leq2$, then $T=\infty$, which implies that $\tilde\mu_0=0$. In this case $x_0=\lambda(\beta-1)$.
\item If $\beta>2$, then $T=(\beta-1)/(\beta-2)$ and, by \eqref{condBAB}, the absorption rate $\tilde\mu_0$ is subject to the restriction
\begin{equation}\label{weywey}
0\leq\tilde\mu_0\leq\lambda(\beta-2).
\end{equation}
\end{itemize}

Under the conditions considered above, the Darboux transformation $\widetilde{\mathcal{A}}$ becomes the infinitesimal generator of a new birth--death process with birth--death rates $\{\tilde{\lambda}_n, \tilde{\mu}_n\}$ given by \eqref{newbdcoe2} and \eqref{laconcha3r}. The corresponding spectral measure of $\widetilde{\mathcal{A}}$ is described by the Geronimus transform of $d\psi$ in \eqref{spmec2e}. For that, we compute the moment $m_{-1}$ of the spectral measure $d\psi$ using formula \eqref{mmm1}, yielding $m_{-1}=A=\frac{1}{\lambda(\beta-1)}$. Hence, the transformed spectral measure of the Darboux transformation takes the form
\begin{equation*}
d\widetilde\psi(x)=\frac{x_0+\tilde\mu_0}{\Gamma(\beta)\lambda^\beta}x^{\beta-2}e^{-x/\lambda}dx+\left(1-\frac{x_0+\tilde\mu_0}{\lambda(\beta-1)}\right)\delta_0(dx).
\end{equation*}
The corresponding orthogonal polynomials can be obtained from \eqref{QDTUL} (see also \eqref{QQbb1}). After some computations we get $\tilde Q_0=1$ and
$$
\tilde Q_n(x)=\frac{s_n}{u_{n-1}\pi_{n-1}}\left[\frac{x_0+\tilde\mu_0}{\lambda(\beta-1)\pi_n}L_n^{(\beta-2)}\left(\frac{x}{\lambda}\right)-\left(1-\frac{x_0+\tilde\mu_0}{\lambda(\beta-1)}\right)\frac{x}{\lambda(\beta+n-1)}L_{n-1}^{(\beta)}\left(\frac{x}{\lambda}\right)\right],\quad n\geq1.
$$
Observe that when $x_0+\tilde\mu_0=\lambda(\beta-1)$ the transformed measure is the Laguerre measure with shape parameter \(\beta-1\), and the corresponding orthogonal polynomials are Laguerre polynomials with parameter \(\beta-2\). In this situation, we can derive an explicit form of the coefficients $x_n, y_n, s_n,$ and $r_n$ in \eqref{xxnn2} and \eqref{yynn2}--\eqref{rrnn2}. Indeed, since \(u_n=\pi_n^{-1}\) in this case, the coefficients take the following forms. For \(1<\beta\leq2\), we necessarily have $\tilde\mu_0=0$, and therefore
\begin{align*}
x_n&=\lambda(\beta-1),\quad n\geq0,\\[0.5em]
y_n&=-x_{n-1},\quad
s_n=\frac{\lambda(n+\beta-1)}{x_{n-1}},\quad
r_n=-\frac{\lambda n}{x_{n-1}},\quad n\geq1.
\end{align*}
For $\beta>2$, we obtain instead
\begin{align*}
x_n
&=\lambda(\beta-1)
-\frac{\tilde\mu_0}{\beta-2}
\left[\beta-1-\frac{(n+1)!}{(\beta)_n}\right],
\quad n\geq0,\\[0.5em]
y_n&=-x_{n-1},\quad
s_n=\frac{\lambda(n+\beta-1)}{x_{n-1}},\quad
r_n=-\frac{\lambda n}{x_{n-1}},\quad n\geq1,
\end{align*}
with \(y_0=-\lambda(\beta-1)\) in both cases. The birth--death rates of the Darboux-transformed process in \eqref{newbdcoe2} and \eqref{laconcha3r} are obtained by treating the boundary separately. Since $x_0+\tilde\mu_0=\lambda(\beta-1)$, we have
\[
\tilde\lambda_0=x_0=\lambda(\beta-1)-\tilde\mu_0,
\quad
\tilde\mu_1=\frac{\lambda^2(\beta-1)}{x_0}.
\]
For $n\geq1$, the remaining rates are
\begin{equation*}
\tilde\lambda_n
=\lambda(n+\beta-1)\frac{x_n}{x_{n-1}},
\quad
\tilde\mu_{n+1}
=\lambda(n+1)\frac{x_{n-1}}{x_n}.
\end{equation*}
If $\tilde\mu_0=0$, then $x_n=\lambda(\beta-1)$ for all $n\geq0$, and therefore
\[
\tilde\lambda_n=\lambda(n+\beta-1),\quad n\geq0,
\quad
\tilde\mu_n=\lambda n,\quad n\geq1.
\]
Thus, the transformed process is precisely the Laguerre birth--death process obtained by replacing $\beta$ with $\beta-1$. On the other hand, when $\beta>2$ and $\tilde\mu_0=\lambda(\beta-2)$ (see \eqref{weywey}), we have $x_n=\lambda(n+1)!/(\beta)_n$. Consequently,
\[
\tilde\lambda_n=\lambda(n+1),\quad n\geq0,
\quad
\tilde\mu_{n+1}=\lambda(n+\beta-1),\quad n\geq0.
\]
This is precisely the dual process (in the sense of Remark~\ref{remdual}) of the Laguerre birth--death process with parameter $\beta-2$; it does not coincide with the dual process of Remark~\ref{remdual2}, since here the original process is conservative.

\section{Linear birth--death processes and associated Laguerre polynomials}\label{secassLag}

In this section we consider the two linear models studied by Ismail, Letessier and Valent in \cite{ILV}. Both have the same birth--death rates
\begin{equation}\label{ratesILV}
\lambda_n=n+\alpha+c+1,\quad n\geq0,\quad \mu_n=n+c,\quad n\geq1,
\end{equation}
and they only differ in the absorption rate at state $0$:
\[
\mbox{Model I:}\quad \mu_0=c,\qquad \mbox{Model II:}\quad \mu_0=0.
\]
We denote by $\mathcal{A}_{\rm I}(\alpha,c)$ and $\mathcal{A}_{\rm II}(\alpha,c)$ the corresponding infinitesimal generators \eqref{QQmm}, and by $d\psi_{\rm I}$ and $d\psi_{\rm II}$ (or $d\psi_{\rm I}^{(\alpha,c)}$ and $d\psi_{\rm II}^{(\alpha,c)}$ when the parameters need to be specified) their spectral measures. To avoid any confusion with the parameter $\alpha$, we will not use here the notation $\mathcal{A}^{(\alpha)}$ of Section~\ref{sec2}; observe however that $\mathcal{A}_{\rm I}(\alpha,c)$ is obtained from $\mathcal{A}_{\rm II}(\alpha,c)$ by replacing $\mu_0=0$ with $\mu_0=c$. 

\begin{figure}[htbp]
\centering
\begin{tikzpicture}[scale=0.9,>=stealth]
\fill[blue!15] (-1,0) -- (-5,4) -- (-5,4.6) -- (3.6,4.6) -- (3.6,0) -- cycle;
\fill[orange!25] (0,-1) -- (-1,0) -- (3.6,0) -- (3.6,-1) -- cycle;
\draw[gray,densely dotted] (-1,0) -- (-1,4.6);
\draw[->] (-5.3,0) -- (4,0) node[right]{$\alpha$};
\draw[->] (0,-1.6) -- (0,4.75) node[left]{$c$};
\foreach \x in {-4,-3,-2,1,2,3} \draw (\x,0.07) -- (\x,-0.07);
\foreach \y in {1,2,3} \draw (0.07,\y) -- (-0.07,\y);
\node[below left,font=\scriptsize] at (-1,0) {$-1$};
\node[below,font=\scriptsize] at (-3,0) {$-3$};
\node[left,font=\scriptsize] at (0,1) {$1$};
\node[left,font=\scriptsize] at (0,-1) {$-1$};
\draw[thick,dashed] (-5,4) -- (0,-1);
\node[font=\footnotesize,rotate=-45] at (-3.35,1.95) {$\alpha+c+1=0$};
\draw[thick,dashed] (0,-1) -- (3.6,-1) node[pos=0.8,below,font=\footnotesize]{$c=-1$};
\draw[very thick,red] (-1,0) -- (3.6,0);
\node[font=\footnotesize,red,above] at (3.1,0) {$c=0$};
\node[font=\footnotesize] at (1.9,2.3) {Models I and II};
\node[font=\footnotesize] at (2.0,-0.55) {Model II only};
\node[font=\scriptsize] at (-3,5.05) {$\alpha<-1$};
\node[font=\scriptsize] at (-0.5,5.05) {$-1\le\alpha\le0$};
\node[font=\scriptsize] at (1.8,5.05) {$\alpha>0$};
\end{tikzpicture}
\caption{Admissible values of the parameters $(\alpha,c)$ in \eqref{ratesILV}.}
\label{fig:ILVparams}
\end{figure}
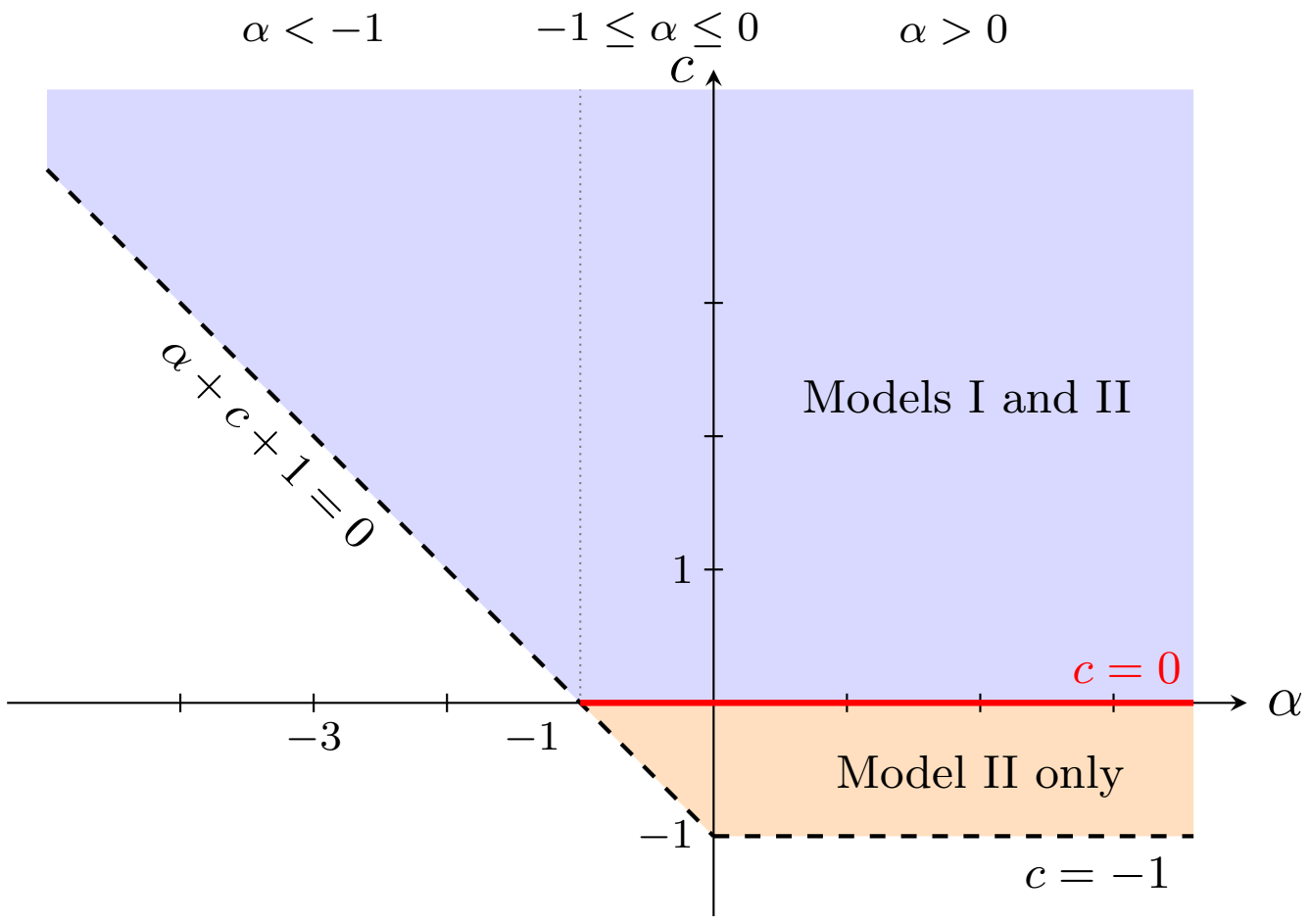

The positivity of the rates requires $\alpha+c+1>0$, together with $c>0$ in Model I and $c>-1$ in Model II (see Figure \ref{fig:ILVparams} for the admissible values of both parameters $\alpha$ and $c$).

If $c=0$, both models coincide with the Laguerre case $\lambda=\mu=1$, $\beta=\alpha+1$, of Section~\ref{seclinear}, so we will assume $c\neq0$. In \cite{ILV} the polynomials $F_n(x)=\pi_nQ_n(x)$ of Model I are the \emph{associated Laguerre polynomials} $L_n^{\alpha}(x;c)$ of Askey and Wimp \cite{AW}, while those of Model II are denoted by $\mathcal{L}_n^{\alpha}(x;c)$.

\subsection{Common quantities}

The potential coefficients \eqref{potcoef} are the same for both models, namely
\begin{equation*}
\pi_n=\frac{(\alpha+c+1)_n}{(c+1)_n},\quad\lambda_n\pi_n=\frac{(\alpha+c+1)_{n+1}}{(c+1)_n},\quad n\geq0,
\end{equation*}
and, as $n\to\infty$, $\pi_n\sim\Gamma(c+1)n^{\alpha}/\Gamma(\alpha+c+1)$. Since $\lambda_n-\mu_n=\alpha+1$ for $n\ge1$, $\lambda_0-c=\alpha+1$, $\lambda_n-\mu_{n+1}=\alpha$ and $\mu_{n+1}-\lambda_{n-1}=1-\alpha$, the symmetry relation $\lambda_{n-1}\pi_{n-1}=\mu_n\pi_n$ gives the telescoping identities
\[
(\alpha+1)\pi_n=\lambda_n\pi_n-\lambda_{n-1}\pi_{n-1},\quad
\frac{\alpha}{\lambda_n\pi_n}=\frac{1}{\pi_n}-\frac{1}{\pi_{n+1}},\quad
\frac{1-\alpha}{\pi_n}=\frac{\mu_{n+1}}{\pi_n}-\frac{\mu_n}{\pi_{n-1}},\quad n\geq0,
\]
(with the conventions $\lambda_{-1}\pi_{-1}=c$ and $\mu_0/\pi_{-1}=\alpha+c$), and therefore
\begin{align}
\label{sumpiILV}\sum_{k=0}^n\pi_k&=\frac{\lambda_n\pi_n-c}{\alpha+1},\quad \alpha\neq-1,\\
\label{sumapiILV}a_n:=\sum_{k=0}^{n-1}\frac{1}{\lambda_k\pi_k}&=\frac{1}{\alpha}\left(1-\frac{1}{\pi_n}\right)=\frac{1}{\alpha}\left(1-\frac{(c+1)_n}{(\alpha+c+1)_n}\right),\quad\alpha\neq0,\\
\label{suminvpiILV}\sum_{k=0}^n\frac{1}{\pi_k}&=\frac{1}{1-\alpha}\left(\frac{\mu_{n+1}}{\pi_n}-\alpha-c\right),\quad\alpha\neq1.
\end{align}
In the exceptional cases $\alpha\in\{-1,0,1\}$, the corresponding sums are expressed in terms of the digamma function $\psi^{(0)}=\Gamma'/\Gamma$: 
\[
\begin{array}{lll}
\alpha=-1: & \pi_k=\dfrac{c}{k+c}, & \D\sum_{k=0}^n\pi_k=c\,\bigl(\psi^{(0)}(n+c+1)-\psi^{(0)}(c)\bigr),\\[4mm]
\alpha=0: & \pi_k=1, & a_n=\psi^{(0)}(n+c+1)-\psi^{(0)}(c+1),\\[3mm]
\alpha=1: & \pi_k=\dfrac{k+c+1}{c+1}, & \D\sum_{k=0}^n\dfrac{1}{\pi_k}=(c+1)\bigl(\psi^{(0)}(n+c+2)-\psi^{(0)}(c+1)\bigr).
\end{array}
\]
Letting $n\to\infty$ in \eqref{sumpiILV} and \eqref{sumapiILV}, and in the first two of these formulas for $\alpha=-1$ and $\alpha=0$ (recall that $\psi^{(0)}(x)\to\infty$ as $x\to\infty$), we obtain that the quantities $A$ and $B$ in \eqref{quantity_AB} are
\begin{equation}\label{ABILV}
A=\begin{cases}
\dfrac{1}{\alpha}, &\mbox{if}\quad \alpha>0,\\[2mm]
\infty, &\mbox{if}\quad \alpha\leq0,
\end{cases}
\qquad
B=\begin{cases}
\dfrac{c}{-\alpha-1}, &\mbox{if}\quad \alpha<-1,\\[2mm]
\infty, &\mbox{if}\quad \alpha\geq-1.
\end{cases}
\end{equation}
Observe that $A$ does not depend on $c$, and that $B<\infty$ forces $c>-\alpha-1>0$. Therefore:
\begin{itemize}[leftmargin=0.25in]
\item Model II is transient if $\alpha>0$, null recurrent if $-1\le\alpha\le0$, and positive recurrent if $\alpha<-1$. In the transient case, by \eqref{sumapiILV}, the probability of ever reaching $0$ is
\begin{equation}\label{hitILV}
\mathbb{P}_n(\tau_0<\infty)=1-\frac{a_n}{A}=\frac{1}{\pi_n}=\frac{(c+1)_n}{(\alpha+c+1)_n},\quad n\ge0.
\end{equation}
In the positive recurrent case, if $X_\pi$ is a random variable with the stationary distribution $\mathbb{P}(X_\pi=n)=\pi_n/B$, then \eqref{sumpiILV} gives the remarkably simple tail
\[
\mathbb{P}(X_\pi>n)=\frac{\lambda_n\pi_n}{c}=\frac{(\alpha+c+1)_{n+1}}{(c)_{n+1}},\quad n\ge0.
\]
Hence $X_\pi$ follows a \emph{Waring (or beta-geometric) distribution}: it is the number of failures before the first success in Bernoulli trials whose success probability is $\mathrm{Beta}(-\alpha-1,\alpha+c+1)$ distributed. In particular, $\mathbb{P}(X_\pi>n)\sim\Gamma(c)\,n^{\alpha+1}/\Gamma(\alpha+c+1)$ as $n\to\infty$, so that, in contrast with the examples of Sections~\ref{secMMinf} and~\ref{seclinear}, the stationary distribution has a power-law tail. Moreover, by \eqref{sumpiILV} and the standard formula $\mathbb{E}_n[\tau_{n+1}]=\frac{1}{\lambda_n\pi_n}\sum_{k=0}^n\pi_k$ (see~\cite[Section 8.1]{And} and~\cite[Section 3.1]{MDIB}), we obtain
\begin{equation}\label{EtauILV}
\mathbb{E}_n[\tau_{n+1}]
=
\begin{cases}
\displaystyle
\frac{1}{\alpha+1}
\left(
1-\frac{(c)_{n+1}}
{(\alpha+c+1)_{n+1}}
\right),
&\mbox{if}\quad \alpha\neq-1,\\[3mm]
\displaystyle
\psi^{(0)}(n+c+1)-\psi^{(0)}(c),
&\mbox{if}\quad \alpha=-1,
\end{cases}
\quad n\ge0.
\end{equation}

\item In Model I, absorption at $-1$ is certain if and only if $\alpha\le0$, and $B<\infty$ if and only if $\alpha<-1$ (since Model~I is nonconservative, it has no stationary distribution, and we avoid the term ergodic here). By \eqref{Qn0}, \eqref{Qn00} and \eqref{sumapiILV}, for $\alpha\neq0$ we have
\begin{equation}\label{Qn0ILV}
Q_n(0)=1+\frac{c}{\alpha}\left(1-\frac{1}{\pi_n}\right)=\frac{(\alpha+c)\pi_n-c}{\alpha\pi_n},\quad Q_n^{(0)}(0)=-\frac{1}{\alpha}\left(1-\frac{1}{\pi_n}\right),\quad n\ge0.
\end{equation}
Hence $\mathbb{P}_i(\tau_n<\tau_{-1})=Q_i(0)/Q_n(0)$, $0\le i\le n$, is explicit, and for $\alpha>0$ the probability of never being absorbed is
\[
\mathbb{P}_i(\tau_{-1}=\infty)=\frac{Q_i(0)}{1+cA}=1-\frac{c}{\alpha+c}\,\frac{(c+1)_i}{(\alpha+c+1)_i},\quad i\ge0.
\]
\end{itemize}
Both models satisfy the standing condition \eqref{uniqueness_assumption}. For Model II this follows from \eqref{ABILV}, since $A<\infty$ requires $\alpha>0$ and $B<\infty$ requires $\alpha<-1$. For Model I, if $\alpha\ge-1$ then $S\ge B=\infty$, while if $\alpha<-1$ then, by \eqref{Qn0ILV}, $\pi_nQ_n^2(0)\sim c^2/(\alpha^2\pi_n)\to\infty$, so that again $S=\infty$. 

As for the spectral measures, it was shown in \cite[Theorem~2.1]{ILV} that their Stieltjes transforms are given, for $p>0$, by
\begin{equation}\label{StieltILV}
\int_0^\infty\frac{d\psi_{\rm I}(x)}{x+p}=\frac{\Psi(c+1,1-\alpha;p)}{\Psi(c,1-\alpha;p)},\qquad
\int_0^\infty\frac{d\psi_{\rm II}(x)}{x+p}=\frac{\Psi(c+1,1-\alpha;p)}{\Psi(c,-\alpha;p)},
\end{equation}
where $\Psi(a,b;z)$ denotes the \emph{confluent hypergeometric function of the second kind} or \emph{Tricomi function}, and that, for $x>0$,
\begin{equation}\label{densILV}
d\psi_{\rm I}(x)=\frac{x^{\alpha}e^{-x}\left|\Psi(c,1-\alpha;xe^{-i\pi})\right|^{-2}}{\Gamma(c+1)\Gamma(\alpha+c+1)}\, dx,\qquad
d\psi_{\rm II}(x)=\frac{x^{\alpha}e^{-x}\left|\Psi(c,-\alpha;xe^{-i\pi})\right|^{-2}}{\Gamma(c+1)\Gamma(\alpha+c+1)}\,dx.
\end{equation}
In \cite[Theorem~2.1]{ILV} both measures are stated to be absolutely continuous for all admissible parameters. This is correct for Model I, but not for Model II in the positive recurrent case $\alpha<-1$, where the spectral measure must have a mass $1/B$ at $x=0$. This situation is admissible in~\cite{ILV} since \(c>-1\) and \(c+\alpha+1>0\) are the only restrictions. Indeed, using $\Psi(a,b;p)\sim\Gamma(b-1)p^{1-b}/\Gamma(a)$ as $p\to0^+$ for $b>1$, \eqref{StieltILV} gives
\[
\lim_{p\to0^+}p\int_0^\infty\frac{d\psi_{\rm II}(x)}{x+p}=\frac{\Gamma(-\alpha)\Gamma(c)}{\Gamma(c+1)\Gamma(-\alpha-1)}=\frac{-\alpha-1}{c}=\frac1B,\quad\alpha<-1,
\]
while the same limit vanishes for $\alpha\geq-1$ and for Model I. Hence
\begin{equation}\label{psiIIcorr}
d\psi_{\rm II}(x)=\frac{x^{\alpha}e^{-x}\left|\Psi(c,-\alpha;xe^{-i\pi})\right|^{-2}}{\Gamma(c+1)\Gamma(\alpha+c+1)}\,\mathbf{1}_{\{x>0\}}\,dx+\frac{-\alpha-1}{c}\,\mathbf{1}_{\{\alpha<-1\}}\,\delta_0(dx).
\end{equation}
This mass will also be recovered below from the Geronimus transformation \eqref{spmec2e} (see \eqref{GerILVmass}). 
Finally, the birth--death polynomials are $Q_n(x)=F_n(x)/\pi_n$, where, for generic parameter values, $F_n(x)$ admits the
hypergeometric representation \cite[(3.4)]{ILV},
\[
F_n(x)=\frac{(\alpha+1)_n}{n!}\sum_{m=0}^n\frac{(-n)_m\,x^m}{(c+1)_m(\alpha+1)_m}\,{}_3F_2\!\left(\!\begin{array}{c}
m-n,\ m+1-\alpha-\eta,\ c \\[-1.5mm]
-\alpha-n,\ c+m+1
\end{array}; 1 \right),\quad n\ge0,
\]
with $\eta=0$ for Model I and $\eta=1$ for Model II, and they satisfy $\int_0^\infty Q_nQ_m\,d\psi=\delta_{nm}/\pi_n$. For Model II the ${}_3F_2$ at $x=0$ is balanced and the Pfaff--Saalsch\"utz formula gives $F_n(0)=\pi_n$, that is, $Q_n(0)=1$, as it should be. The hypergeometric representation above requires
some care at exceptional parameter values. In
particular, when $\alpha$ is a negative integer,
individual Pochhammer symbols and hypergeometric
terms may become singular, even though the
birth--death polynomials remain well defined. At such values the formula is understood through
the corresponding limit in the parameter $\alpha$,
taken for the complete expression after all
necessary cancellations.

\subsection{Model II}
In this case the original process is conservative, so that $Q_n(0)=1$ for all $n\ge0$.

\medskip

\hspace{-.35cm}\underline{\textbf{LU factorization}}. Since $\mu_0=0$, we have $Q_n(0)=1$, $\tilde q_n=-1$ and $\tilde t_n=1-\hat\mu_0/\lambda_0$, $n\ge1$, with $\lambda_0=\alpha+c+1$. By \eqref{sumpiILV} and the corresponding
exceptional formula for $\alpha=-1$,
the numerator in \eqref{ssnn} becomes
\begin{equation}\label{dnLUII}
\tilde d_n:=
\frac{\hat\mu_0}{\lambda_0}
+\left(1-\frac{\hat\mu_0}{\lambda_0}\right)
\sum_{k=0}^n\pi_k
=
\begin{cases}
\displaystyle
\frac{
(\lambda_0-\hat\mu_0)\lambda_n\pi_n
+(\hat\mu_0-c)\lambda_0}
{(\alpha+1)\lambda_0},
&\mbox{if}\quad \alpha\neq-1,\\[3mm]
\displaystyle
\frac{\hat\mu_0}{c}
+(c-\hat\mu_0)
\bigl[\psi^{(0)}(n+c+1)-\psi^{(0)}(c)\bigr],
&\mbox{if}\quad \alpha=-1,
\end{cases}
\end{equation}
so that, by \eqref{ssnn} and \eqref{rrnn}--\eqref{yynn},
\[
\tilde s_n=\frac{\tilde d_n}{\pi_n},\quad\tilde x_n=-\tilde y_n=\frac{\lambda_n\pi_n}{\tilde d_n},\quad n\ge0,\quad \tilde r_n=-\frac{\tilde d_{n-1}}{\pi_n},\quad n\ge1,
\]
and the rates \eqref{newbdcoe} of the Darboux transformation are
\begin{equation}\label{ratesLUII}
\hat\lambda_n=\mu_{n+1}\frac{\tilde d_{n+1}}{\tilde d_n},\quad\hat\mu_{n+1}=\lambda_{n+1}\frac{\tilde d_n}{\tilde d_{n+1}},\quad n\ge0.
\end{equation}
Since $\lambda_0=\alpha+c+1$ and $B=-c/(\alpha+1)$ (see \eqref{ABILV}), Theorem~\ref{theom33} gives the admissible values of $\hat\mu_0$:
\[
0\le\hat\mu_0\le
\begin{cases}
\alpha+c+1, &\mbox{if}\quad \alpha\ge-1,\\[1mm]
c, &\mbox{if}\quad \alpha<-1.
\end{cases}
\]
In all cases $\hat\lambda_n\sim n$, so that the Darboux-transformed process is not explosive, and $\widehat S=\infty$ by \eqref{sww}, so that its spectral measure is uniquely determined. By \eqref{spmec1}, this spectral measure is the Christoffel transform $d\widehat\psi(x)=x\,d\psi_{\rm II}(x)/\lambda_0$, which does not depend on $\hat\mu_0$. The factor $x$ removes the mass that $d\psi_{\rm II}$ may have at $0$ (see \eqref{psiIIcorr}), and comparing the resulting density with \eqref{densILV} we obtain
\begin{equation}\label{ChrILV}
d\widehat\psi(x)=\frac{x}{\alpha+c+1}\,d\psi_{\rm II}^{(\alpha,c)}(x)=d\psi_{\rm I}^{(\alpha+1,c)}(x).
\end{equation}
Thus, for $c>0$, the Christoffel transform of the spectral measure of Model~II with parameters $(\alpha,c)$ is the spectral measure of Model~I with parameters $(\alpha+1,c)$. For $-1<c<0$, Model~I is not defined, and the right-hand side of \eqref{ChrILV} denotes the probability measure whose density is given by the Model~I formula in \eqref{densILV} with parameters $(\alpha+1,c)$.

Since $\tilde y_n=-\tilde x_n$, \eqref{Qnhat} gives $\hat Q_n(x)=\tilde x_n\big(Q_n(x)-Q_{n+1}(x)\big)/x$. Three choices of $\hat\mu_0$ give especially simple results:
\begin{itemize}[leftmargin=0.25in]

\item Case $\hat\mu_0=c$ (admissible if and only if $c\ge0$, i.e. always in the common region of Models I and II; for $\alpha<-1$, where necessarily $c>0$, it is the upper endpoint). Then $\tilde d_n=\lambda_n\pi_n/\lambda_0$ and the factorization is
\[
\tilde s_n=\frac{\lambda_n}{\alpha+c+1},\quad\tilde r_n=-\frac{\mu_n}{\alpha+c+1},\quad\tilde x_n=-\tilde y_n=\alpha+c+1.
\]
By \eqref{ratesLUII}, the Darboux-transformed rates are
\[
\hat\lambda_n=\lambda_{n+1}=n+\alpha+c+2,\quad\hat\mu_n=n+c,\quad n\ge1,\quad\hat\mu_0=c,
\]
that is, $\widehat{\mathcal{A}}=\mathcal{A}_{\rm I}(\alpha+1,c)$, in agreement with \eqref{ChrILV}. As a consequence, the polynomials of the two models are related by
\begin{equation*}
Q_n^{{\rm I},(\alpha+1,c)}(x)=(\alpha+c+1)\,\frac{Q_n^{{\rm II},(\alpha,c)}(x)-Q_{n+1}^{{\rm II},(\alpha,c)}(x)}{x},\quad n\ge0.
\end{equation*}
This explains, in terms of an LU factorization, the connection between the spectral measures of Models I and II observed in \cite{ILV}.

\item Case $\hat\mu_0=\lambda_0=\alpha+c+1$ (always admissible; it is the upper endpoint if $\alpha\ge-1$, and an interior point of $[0,c]$ if $\alpha<-1$, since then $c>\lambda_0$). Then $\tilde d_n=1$, so that $\tilde x_n=-\tilde y_n=\lambda_n\pi_n$, and \eqref{ratesLUII} gives the dual process of Remark~\ref{remdual}, with rates $\hat\lambda_n=\mu_{n+1}=n+c+1$ and $\hat\mu_n=\lambda_n=n+\alpha+c+1$, $n\ge0$, that is, $\widehat{\mathcal{A}}=\mathcal{A}_{\rm I}(-\alpha-1,\alpha+c+1)$. Unlike the case $\hat\mu_0=c$, this is a genuine Model I for all admissible $(\alpha,c)$, since $\alpha+c+1>0$ and $c+1>0$. By \eqref{ChrILV}, Model I with parameters $(\alpha+1,c)$ and with parameters $(-\alpha-1,\alpha+c+1)$ have the same spectral measure, although their rates are different. Indeed, both have $\lambda_n+\mu_n=2n+\alpha+2c+2$ and $\lambda_n\mu_{n+1}=(n+c+1)(n+\alpha+c+2)$, so they are two different splittings of the same Jacobi matrix. This can also be checked from \eqref{densILV} using Kummer's transformation $\Psi(c,-\alpha;z)=z^{\alpha+1}\Psi(\alpha+c+1,\alpha+2;z)$. Accordingly, their polynomials are proportional: $Q_n^{{\rm I},(-\alpha-1,\alpha+c+1)}=\frac{(\alpha+c+2)_n}{(c+1)_n}\,Q_n^{{\rm I},(\alpha+1,c)}$.

\item Case $\hat\mu_0=0$. Then $\tilde d_n=\sum_{k=0}^n\pi_k$, so that, using $\lambda_{n-1}\pi_{n-1}=\mu_n\pi_n$,
\[
\tilde s_n=\lambda_n\mathbb{E}_n[\tau_{n+1}],\quad\tilde r_n=-\mu_n\mathbb{E}_{n-1}[\tau_n],\quad\tilde x_n=-\tilde y_n=\frac{1}{\mathbb{E}_n[\tau_{n+1}]},
\]
with $\mathbb{E}_n[\tau_{n+1}]$ given by \eqref{EtauILV}. The Darboux-transformed rates are, for
$\alpha\neq-1$,
\[
\hat\lambda_n
=\mu_{n+1}
\frac{\lambda_{n+1}\pi_{n+1}-c}
{\lambda_n\pi_n-c},
\quad
\hat\mu_{n+1}
=\lambda_{n+1}
\frac{\lambda_n\pi_n-c}
{\lambda_{n+1}\pi_{n+1}-c},
\quad n\ge0,
\]
which are no longer linear unless $c=0$ (in which case this choice coincides with $\hat\mu_0=c$ and we recover the Laguerre case of Section~\ref{seclinear}). The exceptional case $\alpha=-1$ must be treated
separately. In this case,
\[
\lambda_n=n+c,\quad
\pi_n=\frac{c}{n+c},\quad n\geq0,
\]
and therefore
\(\sum_{k=0}^n\pi_k
=c\bigl[\psi^{(0)}(n+c+1)-\psi^{(0)}(c)\bigr].
\)
Consequently, the Darboux-transformed rates are
$$
\hat\lambda_n
=(n+c+1)
\frac{\psi^{(0)}(n+c+2)-\psi^{(0)}(c)}
{\psi^{(0)}(n+c+1)-\psi^{(0)}(c)},\quad
\hat\mu_{n+1}
=(n+c+1)
\frac{\psi^{(0)}(n+c+1)-\psi^{(0)}(c)}
{\psi^{(0)}(n+c+2)-\psi^{(0)}(c)},
\quad n\ge0.
$$
Thus the LU factorization remains admissible
when $\alpha=-1$, although its rates must be
expressed in terms of the digamma function.
\end{itemize}

\medskip

\hspace{-.35cm}\underline{\textbf{UL factorization}}. Since $\mu_0=0$, condition \eqref{condULgen} requires $A<\infty$, that is, $\alpha>0$, and then $1/A=\alpha$. For $\alpha\le0$ the UL factorization can be applied instead to the dual process, which by the LU factorization above is $\mathcal A_{\rm I}(-\alpha-1,\alpha+c+1)$; this case is covered by the UL factorization of Model~I in Section~5.3. By \eqref{uun}, \eqref{Qn0}, \eqref{Qn00} and \eqref{sumapiILV}, we have
\[
u_n=1-\frac{x_0+\tilde\mu_0}{\alpha}\left(1-\frac{1}{\pi_n}\right)=\left(1-\frac{x_0+\tilde\mu_0}{\alpha}\right)+\frac{x_0+\tilde\mu_0}{\alpha}\,\mathbb{P}_n(\tau_0<\infty),\quad n\ge0,
\]
in agreement with the probabilistic interpretation of $u_n$ given in \cite[Section~4.1]{AdI}, where now $\mathbb{P}_n(\tau_0<\infty)$ is given by \eqref{hitILV}. By \eqref{mmm1} we have $m_{-1}=A=1/\alpha$, so that the Geronimus transform \eqref{spmec2e} is
\begin{equation}\label{GerILVII}
d\widetilde\psi(x)=\frac{x_0+\tilde\mu_0}{x}\,d\psi_{\rm II}(x)+\left(1-\frac{x_0+\tilde\mu_0}{\alpha}\right)\delta_0(dx).
\end{equation}
By Theorem~\ref{theom37}, two cases appear.
\begin{itemize}[leftmargin=0.25in]
\item Case $0<x_0<\alpha$. Then $\tilde\mu_0=0$ and $x_0+\tilde\mu_0=x_0$, and indeed $T=\infty$, since $u_n\to1-x_0/\alpha>0$ and $\pi_n\to\infty$. By \eqref{xxnn2}--\eqref{rrnn2},
\[
x_n=-y_n=\frac{x_0}{\pi_nu_n},\quad n\geq0,\quad s_n=\frac{\lambda_{n-1}\pi_{n-1}u_{n-1}}{x_0},\quad r_n=-\frac{\mu_n\pi_nu_n}{x_0},\quad n\ge1.
\]
The Darboux-transformed rates are $\tilde\lambda_0=x_0$, $\tilde\mu_0=0$ and
\[
\tilde\lambda_n=\mu_n\frac{u_{n-1}}{u_n},\quad n\ge1,\quad\tilde\mu_{n+1}=\lambda_n\frac{u_{n+1}}{u_n},\quad n\ge0,
\]
and the spectral measure is \eqref{GerILVII} with $x_0+\tilde\mu_0=x_0$, which has a positive mass $1-x_0/\alpha$ at $0$. Hence the Darboux-transformed process is positive recurrent, with $\tilde B=\alpha/(\alpha-x_0)$.
\item Case $x_0+\tilde\mu_0=\alpha$. Then $u_n=1/\pi_n=\mathbb{P}_n(\tau_0<\infty)$ and, by Gauss' summation formula, the series \eqref{TTT} is
\begin{equation}\label{TILVII}
T=\sum_{n=0}^\infty\frac{1}{\pi_n}={}_2F_1\!\left(\!\begin{array}{c}
c+1,\ 1 \\[-1.5mm]
\alpha+c+1
\end{array}; 1 \right)=\begin{cases}
\dfrac{\alpha+c}{\alpha-1}, &\mbox{if}\quad \alpha>1,\\[2mm]
\infty, & \mbox{if}\quad0<\alpha\le1,
\end{cases}
\end{equation}
where the divergence for $0<\alpha\le1$ follows from $1/\pi_n\sim\Gamma(\alpha+c+1)n^{-\alpha}/\Gamma(c+1)$. Hence, by \eqref{condBAB}, $\tilde\mu_0=0$ if $0<\alpha\le1$, while for $\alpha>1$ the absorption rate can be chosen in the range
\[
0\le\tilde\mu_0\le\frac{\alpha(\alpha-1)}{\alpha+c}.
\]
In this case $\pi_nu_n=1$, and writing $d_n=\alpha-\tilde\mu_0\sum_{k=0}^n1/\pi_k$ (given explicitly by \eqref{suminvpiILV}) and $d_{-1}=\alpha$, we get
\[
x_n=d_n,\quad y_n=-d_{n-1},\quad n\ge0,\quad s_n=\frac{\lambda_{n-1}}{d_{n-1}},\quad r_n=-\frac{\mu_n}{d_{n-1}},\quad n\ge1,
\]
and the Darboux-transformed rates are $\tilde\lambda_0=x_0$ and
\[
\tilde\lambda_n=\lambda_{n-1}\frac{d_n}{d_{n-1}},\quad\tilde\mu_n=\mu_n\frac{d_{n-2}}{d_{n-1}},\quad n\ge1.
\]
By \eqref{GerILVII}, the spectral measure $d\widetilde\psi=\alpha x^{-1}d\psi_{\rm II}$ has no mass at $0$ and does not depend on $\tilde\mu_0$. Two choices of $\tilde\mu_0$ give especially simple results:
\begin{itemize}[leftmargin=0.25in]
\item If $\tilde\mu_0=0$, then $d_n=\alpha$ and
\[
\tilde\lambda_0=\alpha,\quad\tilde\lambda_n=n+\alpha+c,\quad\tilde\mu_n=n+c,\quad n\ge1,\quad\tilde\mu_0=0,
\]
which is Model II with parameters $(\alpha-1,c)$, except for the birth rate at state $0$, which is $\alpha$ instead of $\alpha+c$ (both coincide when $c=0$, recovering the Laguerre case of Section~\ref{seclinear}).
\item If $\alpha>1$ and $\tilde\mu_0=\alpha(\alpha-1)/(\alpha+c)$, then, by \eqref{suminvpiILV}, $d_n=\alpha\mu_{n+1}/((\alpha+c)\pi_n)$, and
\[
\tilde\lambda_0=\frac{\alpha(c+1)}{\alpha+c},\quad\tilde\lambda_n=n+c+1,\quad\tilde\mu_n=n+\alpha+c-1,\quad n\ge1,\quad\tilde\mu_0=\frac{\alpha(\alpha-1)}{\alpha+c}.
\]
For $n\ge1$ these are the rates \eqref{ratesILV} with parameters $(1-\alpha,\alpha+c-1)$ (which correspond to a genuine Model~I only if $\alpha+c>1$; if $c\le 1-\alpha$, the death rates $\tilde\mu_n=n+\alpha+c-1$ are still positive for $n\ge1$, since $\alpha+c>0$). The boundary rates, however, are not those of Model~I with these parameters, whose absorption rate would be $\alpha+c-1$. Instead, $\tilde\lambda_0$ and $\tilde\mu_0$ are obtained by multiplying the boundary rates $c+1$ and $\alpha-1$ of the analogous case for Model~I below (where $x_0+\tilde\mu_0=\alpha+c$ and $\tilde\mu_0=\alpha-1$) by the common factor $\alpha/(\alpha+c)$, so that $\tilde\lambda_0+\tilde\mu_0=\alpha$, whereas in Model~I one has $\tilde\lambda_0+\tilde\mu_0=\alpha+c$. In other words, the two Darboux-transformed processes have the same rates for $n\ge1$ and differ only at state $0$.
\end{itemize}
\end{itemize}

\subsection{Model I}
In this case the original process is nonconservative, with absorption rate $\mu_0=c>0$.

\medskip

\hspace{-.35cm}\underline{\textbf{LU factorization}}. Now $\lambda_0+\mu_0=\alpha+2c+1$, and $\tilde q_n=-Q_{n+1}(0)/Q_n(0)$, with $Q_n(0)$ given by \eqref{Qn0ILV}. Expanding $\pi_kQ_k^2(0)$ and using \eqref{sumpiILV} and \eqref{suminvpiILV}, the sum $\Sigma_n=\sum_{k=0}^n\pi_kQ_k^2(0)$ appearing in \eqref{ssnn} is, for $\alpha\notin\{0,\pm1\}$,
\begin{equation*}
\Sigma_n=\frac{1}{\alpha^2}\left[\frac{(\alpha+c)^2(\lambda_n\pi_n-c)}{\alpha+1}-2c(\alpha+c)(n+1)+\frac{c^2}{1-\alpha}\left(\frac{\mu_{n+1}}{\pi_n}-\alpha-c\right)\right],\quad n\ge0,
\end{equation*}
The three exceptional cases $\alpha\in\{-1,0,1\}$ can
also be evaluated explicitly. For $\alpha=-1$, we have
\(
\pi_n=c/(n+c)\), \(Q_n(0)=n+1,\) and hence
\begin{equation*}
\Sigma_n=c\left[
(n+1)\left(\frac n2+2-c\right)
+(1-c)^2\bigl(\psi^{(0)}(n+c+1)-\psi^{(0)}(c)\bigr)
\right].
\end{equation*}
For $\alpha=0$, let
\(
h_n=\psi^{(0)}(n+c+1)-\psi^{(0)}(c+1).
\)
Since $\pi_n=1$ and $Q_n(0)=1+ch_n$, we obtain
\begin{equation*}
\Sigma_n={}n+1
+2c\bigl[(n+c+1)h_n-n\bigr]\\
+c^2\bigl[(n+c+1)h_n^2
-(2n+2c+1)h_n+2n\bigr].
\end{equation*}
Finally, for $\alpha=1$,
\[
\pi_n=\frac{n+c+1}{c+1},\quad
Q_n(0)=\frac{(c+1)(n+1)}{n+c+1},\quad n\geq0,
\]
and consequently
\begin{equation*}
\Sigma_n=(c+1)\left[
(n+1)\left(\frac n2+1-c\right)
+c^2\bigl(\psi^{(0)}(n+c+2)-\psi^{(0)}(c+1)\bigr)
\right].
\end{equation*}
Thus $\Sigma_n$ admits an explicit representation
for every admissible pair $(\alpha,c)$.
The coefficients of the factorization follow from \eqref{ssnn} and \eqref{rrnn}--\eqref{yynn}. By Theorem~\ref{theom33}, the admissible range is
\[
0\le\hat\mu_0\le\alpha+2c+1,
\]
and the spectral measure of $\widehat{\mathcal{A}}$ is $d\widehat\psi(x)=x\,d\psi_{\rm I}(x)/(\alpha+2c+1)$. From \eqref{sww} and $\Sigma_n\ge\max\{1,\pi_n Q_n^2(0)\}$, $\widehat S\ge(\lambda_0+\mu_0)\sum_nQ_n(0)/(\lambda_nQ_{n+1}(0))$, which diverges since $Q_n(0)/Q_{n+1}(0)\to1$ and $\sum_n1/\lambda_n=\infty$. This shows that $\widehat S=\infty$, so all these Darboux transforms satisfy \eqref{uniqueness_assumption}. The rates of the Darboux transformations are not linear in general. The simplest case is the endpoint
$\hat\mu_0=\alpha+2c+1$, for which
$\tilde s_n=1/(\pi_nQ_n(0))$ and,
by \eqref{newbdcoe},
\[
\hat\lambda_n
=\mu_{n+1}\frac{Q_n(0)}{Q_{n+1}(0)},
\quad
\hat\mu_n
=\lambda_n\frac{Q_{n+1}(0)}{Q_n(0)},
\quad n\ge0.
\]
For $\alpha\neq0$, these rates can
also be written explicitly as
\[
\hat\lambda_n
=\lambda_n
\frac{(\alpha+c)\pi_n-c}
     {(\alpha+c)\pi_{n+1}-c},
\quad
\hat\mu_n
=\mu_{n+1}
\frac{(\alpha+c)\pi_{n+1}-c}
     {(\alpha+c)\pi_n-c},
\quad n\ge0.
\]
These are precisely the rates of the
dual of the Doob $h$-transform of
Model I with $h(n)=Q_n(0)$. The probabilistic interpretation of
the factorization coefficients given
in \cite[Section~3.1]{AdI}
is now explicit. In particular,
the conditional expected passage times
satisfy
\[
\mathbb E_n[
\tau_{n+1}\mid\tau_{n+1}<\tau_{-1}]
=
\frac{\Sigma_n}
{\lambda_n\pi_nQ_n(0)Q_{n+1}(0)},
\quad n\ge0,
\]
where $\Sigma_n=\sum_{k=0}^n
\pi_kQ_k^2(0)$. Moreover, the expected occupation
times before absorption at $-1$ are
\[
\mathbb E_n[T_n^{(-1)}]
=
\begin{cases}
\displaystyle
\frac{Q_n(0)}{\alpha+c},
& \mbox{if}\quad\alpha>0,\\[2mm]
\displaystyle
\frac{\pi_nQ_n(0)}{c},
& \mbox{if}\quad \alpha\le0.
\end{cases}
\]
For $\alpha\le0$, absorption is
certain, whereas for $\alpha>0$
the process has a positive probability
of escaping to infinity.

\medskip

\hspace{-.35cm}\underline{\textbf{UL factorization}}. By \eqref{uun}, \eqref{Qn0}, \eqref{Qn00}
and \eqref{sumapiILV}, we obtain
\[
u_n=
\begin{cases}
\displaystyle
1+\frac{c-x_0-\tilde\mu_0}{\alpha}
\left(1-\frac{1}{\pi_n}\right),
& \mbox{if}\quad\alpha\neq0,\\[3mm]
\displaystyle
1+(c-x_0-\tilde\mu_0)
\bigl[\psi^{(0)}(n+c+1)-\psi^{(0)}(c+1)\bigr],
&\mbox{if}\quad \alpha=0.
\end{cases}
\]
The exceptional case $\alpha=0$ satisfies $A=B=\infty$,
so Theorem~\ref{theom37} implies that
$\tilde\mu_0=0$ and
\[
0<x_0\le c.
\]
At the endpoint $x_0=c$, we have $u_n=1$ for every
$n\ge0$ and the Darboux transformation coincides
with the dual process $\mathcal A_{\rm II}(-1,c)$. The upper bound in \eqref{condULgen} is $\mu_0+1/A=c+\alpha^+$, where $\alpha^+=\max\{\alpha,0\}$, and by \eqref{mmm1} $m_{-1}=1/(c+\alpha^+)$, so that
\begin{equation}\label{GerILVI}
d\widetilde\psi(x)=\frac{x_0+\tilde\mu_0}{x}\,d\psi_{\rm I}(x)+\left(1-\frac{x_0+\tilde\mu_0}{c+\alpha^+}\right)\delta_0(dx).
\end{equation}
\begin{itemize}[leftmargin=0.25in]
\item Case $0<x_0<c+\alpha^+$. Then $T=\infty$ and $\tilde\mu_0=0$. The coefficients and rates are given by the same formulas as in the first case of the UL factorization of Model II, namely $x_n=-y_n=x_0/(\pi_nu_n)$, $s_n=\lambda_{n-1}\pi_{n-1}u_{n-1}/x_0$, $r_n=-\mu_n\pi_nu_n/x_0$, and $\tilde\lambda_0=x_0$, $\tilde\lambda_n=\mu_nu_{n-1}/u_n$, $\tilde\mu_{n+1}=\lambda_nu_{n+1}/u_n$. Two choices of $x_0$ are especially interesting:
\begin{itemize}[leftmargin=0.2in]
\item If $x_0=c$ (which belongs to this case when $\alpha>0$), then $u_n=1$ and we obtain the dual process of Remark~\ref{remdual2}, with rates $\tilde\lambda_n=\mu_n$, $\tilde\mu_{n+1}=\lambda_n$, that is, $\widetilde{\mathcal{A}}=\mathcal{A}_{\rm II}(-\alpha-1,\alpha+c)$.
\item If $x_0=\alpha+c$ (which belongs to this case when $-c<\alpha<0$), then $u_n=1/\pi_n$ and $x_n=\alpha+c$, so that, as in the first case of the LU factorization of Model II,
\[
\tilde\lambda_n=n+\alpha+c,\quad n\ge0,\quad\tilde\mu_n=n+c,\quad n\ge1,\quad\tilde\mu_0=0,
\]
that is, $\widetilde{\mathcal{A}}=\mathcal{A}_{\rm II}(\alpha-1,c)$. In other words, this UL factorization of $\mathcal{A}_{\rm I}(\alpha,c)$ is the reverse of the LU factorization of $\mathcal{A}_{\rm II}(\alpha-1,c)$ with $\hat\mu_0=c$. By \eqref{GerILVI}, the Geronimus transform recovers the mass at $0$ in \eqref{psiIIcorr}:
\begin{equation}\label{GerILVmass}
d\psi_{\rm II}^{(\alpha-1,c)}(x)=\frac{\alpha+c}{x}\,d\psi_{\rm I}^{(\alpha,c)}(x)+\frac{-\alpha}{c}\,\delta_0(dx),\quad\alpha<0,
\end{equation}
and $-\alpha/c=1/B$ for Model II with parameters $(\alpha-1,c)$.
\end{itemize}
\item Case $x_0+\tilde\mu_0=c+\alpha$ with $\alpha>0$. Then $u_n=1/\pi_n$, and $T$ is again given by \eqref{TILVII}. Hence $\tilde\mu_0=0$ if $0<\alpha\le1$, while for $\alpha>1$ the absorption rate can be chosen in the range
\[
0\le\tilde\mu_0\le\alpha-1.
\]
As in Model II, $\pi_nu_n=1$ and, writing $d_n=\alpha+c-\tilde\mu_0\sum_{k=0}^n1/\pi_k$ and $d_{-1}=\alpha+c$, we have $x_n=d_n$, $y_n=-d_{n-1}$, $s_n=\lambda_{n-1}/d_{n-1}$, $r_n=-\mu_n/d_{n-1}$, and $\tilde\lambda_n=\lambda_{n-1}d_n/d_{n-1}$, $\tilde\mu_n=\mu_nd_{n-2}/d_{n-1}$, $n\ge1$. The spectral measure $d\widetilde\psi=(\alpha+c)x^{-1}d\psi_{\rm I}$ has no mass at $0$. If $\tilde\mu_0=0$, then $x_n=\alpha+c$ and, as above, $\widetilde{\mathcal{A}}=\mathcal{A}_{\rm II}(\alpha-1,c)$, so that \eqref{GerILVmass} also holds for $\alpha>0$ without the mass at $0$. If $\alpha>1$ and $\tilde\mu_0=\alpha-1$, then $d_n=\mu_{n+1}/\pi_n$ by \eqref{suminvpiILV}, and
\[
\tilde\lambda_n=n+c+1,\quad n\ge0,\quad\tilde\mu_n=n+\alpha+c-1,\quad n\ge1,\quad\tilde\mu_0=\alpha-1.
\]
These are the rates \eqref{ratesILV} with parameters $(1-\alpha,\alpha+c-1)$, but with absorption rate $\alpha-1$ instead of $\alpha+c-1$.
\item Case $x_0+\tilde\mu_0=c$ with $\alpha\le0$. Then $u_n=1$ and $T=B$. Hence $\tilde\mu_0=0$ if $-1\le\alpha\le0$ (and we obtain the dual process $\mathcal{A}_{\rm II}(-\alpha-1,\alpha+c)$, as above), while for $\alpha<-1$, by \eqref{ABILV} and \eqref{condBAB}, the absorption rate can be chosen in the range
\[
0\le\tilde\mu_0\le-\alpha-1.
\]
Writing $d_n=c-\tilde\mu_0\sum_{k=0}^n\pi_k$ (given explicitly by \eqref{sumpiILV}) and $d_{-1}=c$, we have $x_n=d_n/\pi_n$, $y_n=-d_{n-1}/\pi_n$, $s_n=\lambda_{n-1}\pi_{n-1}/d_{n-1}$, $r_n=-\mu_n\pi_n/d_{n-1}$, and $\tilde\lambda_n=\mu_nd_n/d_{n-1}$, $\tilde\mu_n=\lambda_{n-1}d_{n-2}/d_{n-1}$, $n\ge1$. The spectral measure $d\widetilde\psi=c\,x^{-1}d\psi_{\rm I}$ has no mass at $0$. If $\tilde\mu_0=-\alpha-1$, then $d_n=\lambda_n\pi_n$ by \eqref{sumpiILV}, so that $x_n=\lambda_n$, and
\[
\tilde\lambda_n=n+\alpha+c+1,\quad n\ge0,\quad\tilde\mu_n=n+c-1,\quad n\ge1,\quad\tilde\mu_0=-\alpha-1,
\]
which are the rates \eqref{ratesILV} with parameters $(\alpha+1,c-1)$, but with absorption rate $-\alpha-1$ instead of $c-1$.
\end{itemize}

\begin{remark}
The previous results show that LU and UL factorizations provide a systematic way of moving between the two models of \cite{ILV}. For instance, if $c>0$,
\[
\mathcal{A}_{\rm II}(\alpha,c)\ \xrightarrow{\ \mathrm{LU},\ \hat\mu_0=c\ }\ \mathcal{A}_{\rm I}(\alpha+1,c)\ \xrightarrow{\ \mathrm{UL},\ x_0=\alpha+c+1,\ \tilde\mu_0=0\ }\ \mathcal{A}_{\rm II}(\alpha,c),
\]
and the same holds with the dual processes, since the LU factorization of $\mathcal{A}_{\rm II}(\alpha,c)$ with $\hat\mu_0=\alpha+c+1$ gives $\mathcal{A}_{\rm I}(-\alpha-1,\alpha+c+1)$, whose UL factorization with $x_0=\alpha+c+1$ and $\tilde\mu_0=0$ gives back $\mathcal{A}_{\rm II}(\alpha,c)$. At the level of spectral measures, these relations are the Christoffel and Geronimus transformations \eqref{ChrILV} and \eqref{GerILVmass}. In particular, the mass at $0$ of $d\psi_{\rm II}$ in the positive recurrent case, which is missing in \cite[Theorem~2.1]{ILV}, appears naturally from the Geronimus transformation.
\end{remark}

\begin{remark}
In \cite[Section~4]{ILV} the authors also consider the associated Meixner and associated Charlier polynomials, corresponding to the rates $\lambda_n=\theta(n+\gamma+\beta)$, $\mu_n=n+\gamma$, and $\lambda_n=a$, $\mu_n=n+\gamma$, respectively, again with $\mu_0=\gamma$ or $\mu_0=0$. In these cases the spectral measures are discrete, and their supports are only given implicitly as the zeros of certain transcendental functions (see \cite[(4.7)--(4.8)]{ILV}). The results of Section~\ref{sec2} still apply, and the spectral measures of the Darboux-transformed processes are the corresponding Christoffel and Geronimus transforms, but the explicit formulas of this section are no longer available. Indeed, they rely on the fact that $\lambda_n-\mu_n$ is constant for $n\ge1$, which makes the partial sums of $\pi_n$ and $1/(\lambda_n\pi_n)$ telescoping, and this fails for the associated Meixner (unless $\theta=1$) and associated Charlier rates.
\end{remark}

\section{Concluding remarks}\label{secConc2}

In this paper we have applied the general results of \cite{AdI} on LU and UL factorizations of birth--death generators to several classical families of birth--death processes with linear rates and explicitly known spectral measures: the $M/M/\infty$ queue (Charlier polynomials), the linear birth--death processes (Meixner and Laguerre polynomials) and the two linear models of \cite{ILV} (associated Laguerre polynomials). In all these cases the coefficients of the factors, the admissible ranges of the parameters and the rates of the Darboux-transformed processes can be written explicitly, and their spectral measures are obtained from the original ones by Christoffel and Geronimus transformations.

\smallskip

Several common features have emerged. First, the Christoffel transformation associated with the LU factorization typically shifts the parameters of the classical families, and different admissible values of $\hat\mu_0$ give processes with different rates but the same spectral measure. In the associated Laguerre case it relates the conservative and the nonconservative models of \cite{ILV}, and the corresponding UL factorizations provide the reverse relation. Second, the Geronimus transformation associated with the UL factorization naturally produces masses at the origin; in particular, it shows that the spectral measure of the conservative model of \cite{ILV} has a mass $1/B$ at $x=0$ in the positive recurrent case, which is missing in \cite[Theorem~2.1]{ILV}. Third, the series $T$ which governs the admissible absorption rates in the UL case may be finite or infinite depending on the parameters, and its behavior reflects the decay of the hitting probabilities $\mathbb{P}_n(\tau_0<\infty)$; for linear rates it may diverge in part of the transient region. Finally, in all the examples the Darboux-transformed processes again have a uniquely determined spectral measure and are not explosive.

\smallskip

Several questions remain open. It would be interesting to study the LU and UL factorizations of the associated Meixner and Charlier models of \cite{ILV}, whose spectral measures are discrete with implicitly defined supports. Another natural class of examples is given by birth--death processes with quadratic rates, such as the models of \cite{ILV2} related to the associated continuous dual Hahn polynomials, for which $\sum_n1/\lambda_n<\infty$ and the determinacy and nonexplosion assumptions of Section~\ref{sec2} are no longer automatic. Cubic and quartic rates, for which explosion may occur and the role of the minimal process becomes essential, as well as $q$-analogues of the families considered here, are also of interest. Finally, iterated LU and UL factorizations, leading to Christoffel and Geronimus transformations of higher order, could be used to construct further families of birth--death processes with explicitly known spectral measures.

\section*{Statement on the use of AI tools}
The authors used AI-based tools to assist in checking computations and consistency, and in editing the text and the \LaTeX{} source. All content was reviewed and verified by the authors, who take full responsibility for it.

\end{document}